 \documentclass[12pt]{amsart}
\font\emailfont=cmtt10

\usepackage{amsmath,amsthm,amsfonts,amscd,flafter,epsf}

\newcommand{\Id}{{{\mathchoice {\rm 1\mskip-4mu l} {\rm 1\mskip-4mu l}
{\rm 1\mskip-4.5mu l} {\rm 1\mskip-5mu l}}}}
\newcommand\commentable[1]{#1}
\newcounter{bean}

\newtheorem{conjecture}{Conjecture}[section]

\newtheorem{theorem}{Theorem}[section]
\newtheorem{prop}[theorem]{Proposition}
\newtheorem{cor}[theorem]{Corollary}
\newtheorem{lemma}[theorem]{Lemma}

\newtheorem{defn}[theorem]{Definition}

\newtheorem{remark}[theorem]{Remark}

\def\endproof{\relax\ifmmode\expandafter\endproofmath\else
  \unskip\nobreak\hfil\penalty50\hskip.75em\hbox{}\nobreak\hfil\bull
  {\parfillskip=0pt \finalhyphendemerits=0 \bigbreak}\fi}
\def\endproofmath$${\eqno\bull$$\bigbreak}
\def\bull{\vbox{\hrule\hbox{\vrule\kern3pt\vbox{\kern6pt}\kern3pt\vrule}\hrule}}

\newcommand{\Q}{\mathbb{Q}}
\newcommand{\R}{\mathbb{R}}

\newcommand{\C}{\mathbb{C}}

\newcommand{\Z}{\mathbb{Z}}

\newcommand{\OneHalf}{\frac{1}{2}}

\newcommand{\CP}[1]{{\mathbb{CP}}^{#1}}

\newcommand{\Zmod}[1]{\Z/{#1}\Z}

\newcommand{\Tr}{\mathrm{Tr}}
\newcommand{\Ker}{\mathrm{Ker}}

\newcommand{\Coker}{\mathrm{Coker}}
\newcommand{\ind}{\mathrm{ind}}
\newcommand{\Ind}{\ind}
\newcommand{\Image}{\mathrm{Im}}
\newcommand{\Real}{\mathrm{Re}}

\newcommand{\Spec}{\mathrm{Spec}}

\newcommand{\cm}{\cdot}

\newcommand{\Nbd}[1]{{\mathrm{nd}}(#1)}
\newcommand{\nbd}[1]{\Nbd{#1}}
\newcommand{\CDisk}{D}

\newcommand{\ModSWfour}{\mathcal{M}}
\newcommand{\ModFlow}{\ModSWfour}
\newcommand{\SW}{SW}

\newcommand{\SpinBunTwo}{S}
\newcommand{\SpinBunTwoP}{\SpinBunTwo^+}
\newcommand{\SpinBunTwoM}{\SpinBunTwo^-}

\newcommand{\Dirac}{\mbox{$\not\!\!D$}}

\newcommand{\Maps}{\mathrm{Map}}

\newcommand{\SpinC}{{\mathrm{Spin}}^c}

\newcommand{\Proj}{\Pi}

\newcommand{\DDt}{\frac{\partial}{\partial t}}
\newcommand{\DDtSq}{\frac{\partial^2}{\partial t^2}}

\newcommand{\goesto}{\mapsto}

\newcommand{\DBar}{\overline{\partial}}

\newcommand{\SpecFlow}{{\mathrm{SF}}}
\newcommand{\SF}{\SpecFlow}

\newcommand\sgn{\sigma}

\newcommand\Wedge{\Lambda}

\newcommand\abuts\Rightarrow
\newcommand\Sym{\mathrm{Sym}}

\newcommand{\Hol}{\mathrm{Hol}}

\newcommand\PD{\mathrm{PD}}
\newcommand\cyl{{cyl}}

\newcommand\Torus[1]{\mathbb{T}^{#1}}

\newcommand\sign{\mathrm{sign}}

\newcommand\Lag{L}

\newcommand\Inv{\theta}
\newcommand\NInv{\widehat\Inv} 
\newcommand\Ntheta{\NInv}
\newcommand\Corr{\xi}
\newcommand\CorrY{\Corr}
\newcommand\OurJac{\widetilde{\Jac}}
\newcommand\OurSym{{\widetilde\Sym}^{g-1}(\Sigma)}

\newcommand\OurTheta{\widetilde\Theta}

\newcommand\spinc{\mathfrak a}
\newcommand\spinca\spinc
\newcommand\spincb{\mathfrak b}
\newcommand\spincc{\mathfrak c}
\newcommand\spincd{\mathfrak d}
\newcommand\spincr{\mathfrak x}

\newcommand\Jac{J}

\newcommand\NumDisks{n}

\newcommand{\Tors}{\mathrm{Tors}}
\newcommand{\Umet}{k}
\newcommand{\metU}{\Umet}

\newcommand{\PCMet}{\Umet_0^+}
\newcommand{\PCA}{A_0^+}
\newcommand{\Vmet}{\mu}
\newcommand{\met}{h}
\newcommand{\Ymet}{k_Y}
\newcommand{\metY}{\Ymet}
\newcommand{\Kmet}{k'}
\newcommand{\metC}[1]{k_{#1}^\cyl}
\newcommand{\spincX}{\mathfrak r}

\newcommand{\Harm}{\mathcal H}

\newcommand{\Met}{\mathrm{Met}}

\newcommand{\ModSp}{M}

\newcommand{\Conf}{\mathcal C}
\newcommand{\ConfU}{\Conf^{\circ}}
\newcommand{\CorrU}{\Corr^{\circ}}

\newcommand{\Restrict}{\rho}

\newcommand{\DiracTwo}{\mbox{$\not\!{\partial}$}}
\newcommand{\DiracThree}{\mbox{$\not\!\!{D}$}}

\newcommand\scale{\mu}

\newcommand{\TConn}{\mbox{$^\circ\nabla$}}
\newcommand{\LConn}{\mbox{$\nabla$}}
\newcommand{\PDiracThree}{\mbox{$^\circ\DiracThree$}}
\newcommand\Cyl{\R\times\Sigma}
\newcommand\HDisk{B^+}
\newcommand\TBall{B}
\newcommand\Surf{F}

\newcommand\data{\wp}
\newcommand\PCdata{\data_0^+}

\newcommand\etaDirac{\eta^{\mathrm{Dirac}}}
\newcommand\etaSign{\eta^{sign}}
\newcommand\tSpinC{\underline{\SpinC}}
\newcommand\Dom{C}

\newcommand\ConstOne{L}
\newcommand\ConstTwo{M}

\commentable{
\title{The Theta Divisor and the Casson-Walker Invariant}
\author[Peter Ozsv{\'a}th]{Peter Ozsv\'ath}
\thanks{The first author was partially supported by NSF grant number
9971950}
\address{Department of
Mathematics,  Princeton University, New Jersey 08540 \newline
\indent{\emailfont{petero@math.princeton.edu}}}

\author[Zolt{\'a}n Szab{\'o}]{Zolt{\'a}n Szab{\'o}} 
\thanks{The second author was partially 
supported by NSF grant number DMS 970435, a
Sloan Fellowship, and a Packard Foundation grant} 
\address{Department of
Mathematics,  Princeton University, New Jersey 08540 \newline
\indent{\emailfont{szabo@math.princeton.edu}}}}

\begin{document}

\begin{abstract}
    We use Heegaard decompositions and the theta divisor on
    a Riemannian surface to define a three-manifold invariant for
    rational homology three-spheres.  This invariant is defined on the set
    of $\SpinC$ structures $$\Ntheta\colon \SpinC(Y)\longrightarrow \Q.$$ In
    the first part of the paper, we give the definition of the invariant
    (which builds on the theory developed in~\cite{Theta}). In the second
    part, we establish a relationship between this invariant and the
    Casson-Walker invariant.
\end{abstract}


\maketitle

\section{Introduction}
\label{sec:Introduction}

In~\cite{Theta}, we studied a topological invariant associated to
Heegaard decompositions of oriented three-manifolds whose first Betti
number $b_1(Y)$ is positive. Starting with a Heegaard decomposition
$U_0\cup_{\Sigma}U_1$ for $Y$, the invariant measures how the theta
divisor of $\Sigma$ moves as the surface undergoes degenerations
naturally associated to the Heegaard splitting.  In this paper, we
describe a related construction which works in the case where
$b_1(Y)=0$, giving a function $$\NInv\colon
\SpinC(Y)\longrightarrow \Q$$ 
on the set of $\SpinC$ structures on $Y$. Once again, the invariant
measures how the theta divisor moves under the prescribed
degenerations, except that now a subtlety arises due to a
path-dependence of the earlier construction.

To recall the definition of the invariant from~\cite{Theta}, we give
some geometric background. A Heegaard decomposition $U_0\cup_{\Sigma}
U_1$ of $Y$ is a decomposition of the three-manifold as a union of two
handlebodies, identified along an oriented surface $\Sigma$ of genus
$g$. A handlebody $U$ can be described by attaching to $\Sigma$ $g$
two-handles and one three-handle. Let $\{\gamma_1,...,\gamma_g\}$
denote the attaching circles for these two-handles. As explained
in~\cite{Theta}, the handlebody gives rise to a special class of
metrics on $\Sigma$, the {\em $U$-allowable metrics}. Informally, any
metric which is sufficiently stretched out normal to
$\{\gamma_1,...,\gamma_g\}$ is a $U$-allowable metric. We think of the
Jacobian, $\Jac$, as the space of holomorphic line bundles over
$\Sigma$ with degree $g-1$, which in turn is identified with
$H^1(\Sigma;S^1)$, by specifying a spin structure on $\Sigma$. The
handlebody $U$ gives rise to a canonical $g$-dimensional torus
$L(U)\subset \Jac$, which corresponds to the torus $H^1(U;S^1)\subset
H^1(\Sigma;S^1)$ under the  correspondence induced from any spin
structure on $\Sigma$ which extends to $U$. (Note that $L(U)$ is
independent of the spin structure on $U$.)  A metric $\met$ on $\Sigma$
gives rise to an Abel-Jacobi map $$\Theta\colon
\Sym^{g-1}(\Sigma)\longrightarrow
\Jac,$$ which associates 
to a point $D$ in the $(g-1)$-fold symmetric product the holomorphic
line bundle which admits a holomorphic section vanishing exactly at
$D$.  For a $U$-allowable metric, $L(U)$ is disjoint from
the theta divisor.  Fix a path $\met_t$ such
that $\met_0$ is $U_0$-allowable and $\met_1$ is $U_1$-allowable, and
consider the moduli space $$
\ModSp(\met_t)=\{(s,t,D)\in[0,1]\times[0,1]\times
\Sym^{g-1}(\Sigma)\big| s\leq t,
\Theta_{\met_s}(D)\in \Lag(U_0), \Theta_{\met_t}\in \Lag(U_1)\}.$$ For a
generic path (and small perturbations of the $\Lag(U_0)$ and
$\Lag(U_1)$), these points form a discrete set, which misses the locus
where $s=t$. Moreover, the points can be naturally partitioned
according to $\SpinC$ structures over $Y$, and we define
$\theta_{\met_t}(\spinc)$ to be the signed number of points in the
subset corresponding to the $\SpinC$ structure $\spinc$. When $b_1(Y)>0$,
this signed count is shown to be independent of the path of metrics
and, indeed, independent of the Heegaard decomposition used in its
definition.

By contrast, in the case when $b_1(Y)=0$, the signed count
$\theta_{\met_t}$ depends on the choice of path $\met_t$, since the
moduli space can hit the locus where $s=t$ in a one-parameter family
of paths.  In order to get a well-defined topological invariant, we
correct by certain rational-valued, metric-dependent correction terms.
Specifically, the path of metrics $\met_t$ can be used to construct a
metric on the cylinder $\R\times\Sigma$ (where we stretch out the time
directions, in a manner made precise in Section~\ref{sec:ALimit}), and
the path dependence of the invariant $\theta_{\met_t}(\spinc)$ can be
reinterpreted in terms of spectral flow: if $\met_t$ and $\met_t'$ are
two paths of metrics which connect $U_0$-allowable to $U_1$-allowable
metrics, then the difference in the $\theta$ is related to the 
(complex) spectral flow of the $\SpinC$ Dirac operator over the cylinder
$\R\times\Sigma$ by:
$$\theta_{\met_t'}(\spinc)-\theta_{\met_t}(\spinc)=
\SF(\met_t,\met_t')$$
(see Propositions~\ref{prop:ThetaChanges} and \ref{prop:SFCyl}).  Over
a handlebody $U$, in Section~\ref{sec:Handlebodies}, we describe a
canonical (integer-valued) correction term $\CorrU$ on the space of
metrics over $U$ with $U$-allowable boundary, which also changes by
spectral flow.  Thus, if we take a path of metrics $\met_t$, extend it
over $U_0$ by $\Umet_0$, and $U_1$ by $\Umet_1$, we obtain a metric
$\metY$ over $Y$; and the quantity
$$\CorrU(\Umet_0)+\theta_{\met_t}(\spinc)+\CorrU(\Umet_1)$$ will
depend on the metric $\metY$ only through the spectral flow of its
associated Dirac operator, according to standard splitting results for
the spectral flow. To get a topological invariant, then, it suffices
to subtract off any other quantity which depends on the metric only
through spectral flow in the same manner.

Such a canonical metric-dependent term $\Corr_{\metY}(\spinc)$ is
furnished by the index theory for manifolds with cylindrical ends,
developed by Atiyah-Patodi-Singer~\cite{APSI}. It is defined as follows. 
Choose any four-manifold $X$
bounding $Y$, equipped with a cylindrical-end metric $g_X$ and a
$\SpinC$ structure $\spincX$ which bounds $\Umet_Y$ and $\spinc$
respectively.  Then,
$$\Corr_{\metY}(\spinc)=\ind_{X}(\Dirac)-\left(\frac{c_1(\spincX)^2-\sgn(X)}{8}\right),$$
where $\ind_X(\Dirac)$ is the complex index of the Dirac operator
associated to the $\SpinC$ structure $\spincX$, and $\sgn(X)$ is the
signature of the intersection form of $X$.
By the Atiyah-Singer index theorem, the correction term is independent
of the choice of four-manifold $X$ and extending $\SpinC$ structure
$\spincX$.  Note that $\Corr_{\metY}(\spinc)$ is {\em a priori} a
rational number (since $c_1(\spincX)^2$ is rational). (Equivalently,
$\Corr_{\metY}(\spinc)$ can be obtained as a combination of APS
eta-functions for the Dirac and signature operators -- see
Equation~\eqref{eq:ReexpressCorrY}. This latter point of view is 
exploited in Section~\ref{sec:Lens}.)

We now define the normalized invariant $\NInv$ by the equation
$$\NInv(\spinc)=\CorrU(\Umet_0)+\theta_{\met_t}(\spinc)+\CorrU(\Umet_1)-
\Corr_{\metY}(\spinc).$$ This
construction parallels, and was motivated by, a similar treatment of
the Seiberg-Witten invariant for homology three-spheres
(see~\cite{DonaldsonSW}, \cite{Lim}, \cite{Chen}).

Our first result, whose proof occupies
Sections~\ref{sec:DefNTheta}-\ref{sec:GenMet}, is that the quantity
$\NInv$, whose definition involves certain choices (a Heegaard
decomposition, a family of metrics, etc.), gives a well-defined
three-manifold invariant:

\begin{theorem}
\label{thm:WellDefined}
The function
$$\NInv\colon
\SpinC(Y)\longrightarrow \Q$$ 
is a well-defined topological invariant; in particular, it does not
depend on the metrics, Heegaard decompositions of $Y$.
\end{theorem}

Moreover, we will work out a surgery formula for the invariant, which
gives a relationship between $\NInv$ and the invariant $\theta$ for
manifolds with $b_1(Y)=1$. The surgery formula, and our previous
computations of $\theta$ when $b_1(Y)=1$ from~\cite{Theta}, give the
following link between the complex geometry of the Heegaard
decomposition and the $SU(2)$ representations of the fundamental group
of $Y$:

\begin{theorem}
\label{thm:Casson}
Let $Y$ be an integral homology three-sphere. Then,
$$2~\NInv(Y)=\lambda(Y),$$ where $\NInv(Y)$ is the invariant evaluated
on the unique $\SpinC$ structure of $Y$, and $\lambda(Y)$ is Casson's
invariant normalized so that $4\lambda(Y)\equiv \sign(X) \pmod{16}$
for each spin four-manifold $X$ which bounds $Y$.
\end{theorem}

The surgery formula for integral homology three-spheres, and the proof
of the above theorem, are given in Section~\ref{sec:Casson}.  In
Section~\ref{sec:Walker}, we study the $\NInv$-invariant for rational
homology three-spheres. The main result of that section gives a
relationship between $\NInv$ and Walker's generalization of Casson's
invariant (see~\cite{Walker}):

\begin{theorem}
\label{thm:ThetaWalker}
Let $Y$ be a rational homology three-sphere. Then, 
$$2\sum_{\spinc\in\SpinC(Y)}\NInv(\spinc)=\big|H_1(Y;\Z)\big|\lambda(Y),$$
where $\lambda(Y)$ is the Casson-Walker invariant of $Y$.
\end{theorem}

%

It is interesting to compare the above with the situation in gauge
theory.  Counting solutions to the Seiberg-Witten equation, one
obtains a metric-dependent quantity for rational homology spheres,
which depends on its metric through the spectral flow of the Dirac
operator. This is the signed count $\SW_Y(\spinc)$ of the irreducible
Seiberg-Witten monopoles. This, too, can be corrected by the
metric-dependent quantity $\xi_{\metY}$ to obtain a rational-valued
function $$\SW\colon \SpinC(Y)\longrightarrow \Q.$$ Results of
~\cite{Lim} and
\cite{Chen} show that for integral homology three-spheres 
this invariant agrees with half of Casson's invariant. This, together
with Theorem~\ref{thm:Casson}, underlines the close relationship
between $\NInv$ and the Seiberg-Witten invariant (see
also~\cite{Theta} for the case where $b_1(Y)>0$). 
In fact, we discovered the invariants $\theta$ and $\NInv$ by studying
Seiberg-Witten theory and Heegaard decompositions, and it is very
natural to make the following conjecture:

\begin{conjecture}
The invariant $\NInv$ agrees with the Seiberg-Witten invariant $\SW$
for all rational homology three-spheres.
\end{conjecture}

There are two routes for establishing this conjecture: one is to
compare the surgery formulas for both invariants, another is to
proceed more directly via an adiabatic limit of the Seiberg-Witten
equations. Carrying out either programme would take us rather far from
the scope of the present paper. We hope to return to these topics in a
future paper.


\section{Definition of the Invariant}
\label{sec:DefNTheta}

In this section, we start with studying the metric dependence of the
quantity $\Inv$, with a view to showing the topological invariance of
$\NInv$. We begin with a few preliminary remarks on the
definition of $\Inv$.

\begin{defn}
A path $\{h_t\}_{t\in[0,1]}$ of metrics over $\Sigma$ is called {\em
$(U_0,U_1)$-allowable} if the metric $\met_0$ is $U_0$-allowable, and
the metric $\met_1$ is $U_1$-allowable. 
\end{defn}

Let $\ModSp(\met_t)$ denote the moduli space introduced in
Section~\ref{sec:Introduction}.

\begin{prop}
\label{prop:OurGenMet}
The moduli space $\ModSp(\met_t)$ can be naturally partitioned into
components $\ModSp(\met_t,\spinc)$ labeled by $\SpinC$ structures over
$Y$. Moreover, for a generic $(U_0,U_1)$-allowable path of metrics
$\met_t$ on $\Sigma$, the moduli space is a compact, oriented
$0$-manifold which does not contain any points with $s=t$.
\end{prop}

\begin{proof}
The orientation and partitioning into $\SpinC$ structures are
described in Section~\ref{Theta:sec:DefTheta} of~\cite{Theta}. When
$g>1$, the genericity statement follows from
Proposition~\ref{prop:GenericMetrics}, which is proved in
Section~\ref{sec:GenMet} of the present paper. Note that solutions
with $s=t$ correspond to points in $\Sym^{g-1}(\Sigma)$ which map to
$\Lag(U_0)\cap \Lag(U_1)$, which is a discrete set of points (since
$Y$ is a homology three-sphere); so these are excluded by dimension
counting.

When $g=1$, it is easy to see that the
moduli space is empty: in that case, the theta divisor consists of a
single, isolated spin structure, which does not bound.
\end{proof}

Armed with Proposition~\ref{prop:OurGenMet}, we can define
$\theta_{\met_t}(\spinc)$ to be the signed number of points in the component of
$\ModSp(\met_t,\spinc)$. 

\begin{defn}
Given a pair of
$(U_0,U_1)$-allowable paths $\met_t$ and $\met_t'$, a
$(U_0,U_1)$-allowable homotopy from $\met_t$ to $\met_t'$ is a smooth,
two-parameter family of metrics $$H\colon[0,1]\times[0,1]\rightarrow
\Met(\Sigma)$$ so that for each $t\in[0,1]$, $H(0,t)=\met_t$ and
$H(1,t)=\met_t'$, and for each $u\in[0,1]$ the path $t\mapsto H(u,t)$
is a $(U_0,U_1)$-allowable path.
\end{defn}

We will drop the handlebodies from the notation when they are clear
from the context.  Since the space of $U_i$-allowable metrics is
path-connected for $i=0,1$, see~\cite{Theta} (and the space of metrics
over $\Sigma$ is simply-connected), any two $(U_0,U_1)$-allowable
paths can be connected by a $(U_0,U_1)$-allowable homotopy.

Note that since $Y$ is a rational homology sphere, the $\SpinC$
structures naturally correspond to the intersection points of
$L_0=L(U_0)$ and $L_1=L(U_1)$, see also~\cite{Theta}.

\begin{prop}
\label{prop:ThetaChanges}
Let $\met_t$ and $\met_t'$ be a pair of generic $(U_0,U_1)$-allowable
paths. Then,
$$\theta_{\met_t}(\spinc)-\theta_{\met_t'}(\spinc)=
\#\{(D,u,t)\in\Sym^{g-1}(\Sigma)\times[0,1]\times[0,1]\big|
\Theta_{H(u,t)}(D)=p\},$$
where $H$ is any allowable homotopy from $\met_t$ to $\met_t'$, and
$p\in \Lag_0\cap \Lag_1$ is the point corresponding to the $\SpinC$
structure $\spinc$.
\end{prop}

\begin{remark}
Note that the set
$\{(D,u,t)\in\Sym^{g-1}(\Sigma)\times[0,1]\times[0,1]\big|
\Theta_{H(u,t)}(D)=p\}$ is not necessarily transversally cut out, as
one can easily see by considering spin structures (by Serre duality,
it is easy to see that if $p$ corresponds to a spin structure, then
the local multiplicities are always even). Thus, the intersection
number is to be interpreted by perturbing $p$ slightly.
\end{remark}

\begin{proof}
Let $H$ be a generic allowable homotopy connecting $\met_t$ and
$\met_t'$.  Consider the moduli space $$M(H)= \left\{(D,u,s,t)|s\leq
t, \Theta_{H_{u,s}}(D)\in \Lag(U_0),
\Theta_{H_{u,t}}(D)\in \Lag(U_1)\right\}.$$ According to the generic
metrics result, Proposition~\ref{prop:GenericMetrics}, this space is
an oriented one-manifold with boundary. Since $H$ is an allowable
homotopy, the only boundary components are $u=0$, $u=1$ and
$s=t$. Thus, counting boundaries, with sign and multiplicity, we get
the result as stated.
\end{proof}

The difference term appearing above has another natural interpretation
as a spectral flow on the cylinder $\R\times\Sigma$, inspired by work
of Yoshida (\cite{YoshidaII}, see also~\cite{KK}, \cite{CLMII}).  An
allowable path $\met_t$ and a scale factor $\scale$
naturally induces a metric on the cylinder
$\R\times\Sigma$ given by $(\scale dt)^2 + \met_t$,
where we extend $\met_t$ to all $t\in\R$ by requiring
$\met_t=\met_0$ (resp. $\met_1$) for $t\leq 0$ (resp. $\geq 1$). To
describe the relevant spectral flow, we introduce some 
terminology. 

\begin{defn}
\label{def:Reducible}
A connection $A$ on a spinor bundle $W$ over $Y$ is called {\em
reducible} if the trace of its curvature (or, equivalently, the
curvature induced on its determinant line bundle) vanishes. 
\end{defn}

On a rational homology sphere $Y$ equipped with a Riemannian metric,
each $\SpinC$ structure has a unique reducible connection (up to
gauge). 

\begin{prop}
\label{prop:SFCyl}
Fix a $\SpinC$ structure $\spinc\in \SpinC(Y)$.  Let $\met_t$ be an
allowable path for some Heegaard decomposition of a rational homology
sphere $Y$. Let $A$ be a reducible over $[0,1]\times\Sigma$ which
extends to a reducible over $Y$ for the $\SpinC$ structure $\spinc$.
Then, the Dirac operator coupled to $A$ acting on
$L^2(\R\times\Sigma)$ is Fredholm. If for each $t\in [0,1]$, the theta
divisor for $\met_t$ misses the point corresponding to $\spinc$, then
for all sufficiently large $\scale$, the Dirac operator for the metric
$(\scale dt)^2+\met_t$ and reducible belonging to $\spinc$ has no $L^2$
kernel over $\R\times \Sigma$.  Moreover, let $H$ be an allowable homotopy from
$\met_t$, $\met_t'$ over $\Sigma$. Then, there is some $\scale_0\geq
0$ such that for all sufficiently large $\scale,\scale'\geq \scale_0$,
the spectral flow between the two induced Dirac operators on $\Cyl$
coupled to reducible connections obtained by restricting the
reducibles in $\spinc\in\SpinC(Y)$ is given by
\begin{eqnarray*}
\lefteqn{\SF\left((\mu dt)^2+\met_t, (\mu'dt)^2 + \met_t'\right)} \\
&=& \#\left\{(D,u,t)\in\Sym^{g-1}(\Sigma)\times[0,1]\times[0,1]\big|
\Theta_{H(u,t)}(D)=p\right\},
\end{eqnarray*}
where $p$ is the point in $\Jac(\Sigma)$ corresponding to the $\SpinC$
structure $\spinc$.
\end{prop}

The above is a special case of the more general result
(Proposition~\ref{prop:SFCylBig}), which is proved in
Section~\ref{sec:ALimit}. 

Together, Propositions~\ref{prop:ThetaChanges} and \ref{prop:SFCyl}
show that $\theta_{\met_t}(\spinc)$ depends on the family of metrics
and $\SpinC$ connection used on the cylinder only through the spectral
flow of its Dirac operator. We call metric-dependent quantities which
depend on the metric through only the spectral flow of its associated
Dirac operator {\em chambered metric invariants}.  We use $\theta$,
together with another invariant  $\CorrU$ for handlebodies, to construct a
chambered invariant for certain metrics on $Y$.  Note that an allowable path
$\met_t$ and extensions $\Umet_0$ and $\Umet_1$ of the metrics
$\met_0$ and $\met_1$ over $U_0$ and $U_1$ respectively can be spliced
together naturally to form a metric $(\Umet_0)\#(\met_t)\#(\Umet_1)$
over $Y$.

In Section~\ref{sec:Handlebodies}, we construct a canonical chambered
invariant $\CorrU$ for handlebodies of arbitrary genus, which is uniquely
determined by an excision property, and a normalization condition
which guarantees that $\CorrU(\CDisk^3)=0$ for a non-negative
sectional curvature metric on the three-disk (see
Proposition~\ref{prop:CorrU} for a precise statement of the excision
property). Strictly speaking, this function depends on the metric and
a connection used on the spinor bundle over $U$, and an exact
perturbation of a reducible connection over $U$.  
Thus, the expression $\CorrU(\spinc|U_i,\Umet_i,a_i)$, denotes
the invariant evaluated on the metric $\Umet_i$ and the $\SpinC$
connection obtained by restricting the reducible in $\spinc$ to $U_i$,
and perturbing by $a_i$, where $a_i$ is a one-form which is compactly
supported in the interior of $U_i$. (This one-form is used to ensure
that the kernel of the Dirac operator on the handlebody has no kernel,
see Lemma~\ref{lemma:UDirac}.)

Splitting properties of spectral flow, then, show that the quantity
$\CorrU(\spinc|U_0,\Umet_0,a_0)+\theta_{\met_t}(\spinc)+\CorrU(\spinc|U_1,\Umet_1,a_1)$
depends on the metric $(\Umet_0)\#(\met_t)\#(\Umet_1)$ and the
connections used only through the spectral flow of the associated
Dirac operator, as described in the next proposition. To state the
proposition properly, we must analyze the choices made in ``splicing''
metrics to obtain a metric on $Y$.

Given metrics $\Umet_0$, $\Umet_1$ on the handlebodies $U_0$ and
$U_1$, and a path $\met_t$ on $\Sigma$ which interpolates between the
metrics on $\Sigma$ induced on the boundaries of $U_0$ and $U_1$, we
can splice to obtain a metric on $Y$. The spliced metric requires
three additional parameters: a scale $\scale$ to be used on the path
of metrics $\met_t$, and a pair of ``neck-length'' parameters which
specify lengths of cylinders to be spliced before and after the middle
neck.  More precisely, the metric $$(\Umet_0)\#_{T_0}\left((\scale dt)^2 +
\met_t\right)\#_{T_1} (\Umet_1)$$ will denote the metric on $Y$ which
is obtained by inserting cylinders $[0,T_i]\times \Sigma$ given the
product metric $dt^2 + \met_i$ for $i=0,1$ between the handlebodies
$U_0$ and $U_1$ and the cylinder $[0,1]\times\Sigma$ endowed with the
metric $(\scale dt)^2 + \met_t$.

\begin{prop}
\label{prop:ExtendToY}
Fix a $\SpinC$ structure $\spinc\in\SpinC(Y)$.
Let $\met_t$ be a generic $(U_0,U_1)$-allowable path. For a scale $\scale$ and
necklength parameters $T_0$ and $T_1$, let $\metY(\scale,T_0,T_1)$
denote the metric on $Y$ obtained by splicing 
$$\Umet_0\#_{T_0}\left(\left(\mu dt\right)^2 + \met_t\right) \#_{T_1} \Umet_1.
$$
Let $A(\mu,T_0,T_1)$ denote the corresponding reducible
connections. Then, for generic, compactly supported
one-forms $a_0$, $a_1$ in $U_0$ and $U_1$, the Dirac operator on $Y$
for the metric $\metY(\scale,T_0,T_1)$ and connection
$A(\mu,T_0,T_1)+a_0+a_1$ has no kernel, provided that $\mu$, $T_0$,
and $T_1$ are sufficiently large. Moreover, suppose that $H$ is an allowable
homotopy from $\met_t$ to another allowable path $\met_t'$ which
interpolates between metrics $\Umet_0'$ and $\Umet_1'$. Then, for
generic compactly-supported one-forms  $a_0$, $a_0'$ and $a_1$,
$a_1'$ on $U_0$ and $U_1$ respectively we have that for all 
sufficiently large $\scale$ and necklength parameters $T_0,T_1$ the
spectral flow of the Dirac operator $\Dirac(\scale,T_0,T_1)$
for the metric
$\metY(\scale,T_0,T_1)$ and connection $A(\mu,T_0,T_1)+a_0+a_1$
to the Dirac operator $\Dirac'(\scale,T_0,T_1)$ for the metric 
$\metY'(\scale,T_0,T_1)$ and connection $A'(\mu,T_0,T_1)+a_0'+a_1'$
is given by
\begin{eqnarray*}
\lefteqn{\SF(\Dirac(\scale,T_0,T_1),\Dirac'(\scale,T_0,T_1))} && \\
&=&
\left(\CorrU(\spinc|U_0,\Umet_0',a_0')+
\theta_{\met_t'}(\spinc)+\CorrU(\spinc|U_1,\Umet_1',a_1')\right) \\
&& -
\big(\CorrU(\spinc|U_0,\Umet_0,a_0)+\theta_{\met_t}(\spinc)+
\CorrU(\spinc|U_1,\Umet_1,a_1)\big).
\end{eqnarray*}
\end{prop}

\begin{proof}
The vanishing of the kernel follows from the vanishing of the kernel
over the three pieces: $U_0$, $[0,1]\times \Sigma$, and $U_1$
respectively. Over the handlebodies the kernel vanishes for generic
choices of $a_0$ and $a_1$ (this will be proved in
Lemma~\ref{lemma:UDirac}). Over the cylinder, it vanishes for generic
paths $\met_t$, according to Proposition~\ref{prop:SFCyl}.  Moreover,
the spectral flow statement follows from the splitting principle, the
chambered property of $\CorrU$, and the chambered property of
$\theta$, which in turn follows from
Proposition~\ref{prop:ThetaChanges} together with
Proposition~\ref{prop:SFCyl}.
\end{proof}

In the above proposition, the connection $A=A(\mu,T_0,T_1)+a_0+a_1$
has no kernel for generic $a_0$ and $a_1$. Thus, if we consider a
four-manifold $X$, equipped with a cylindrical-end metric $g_X$ and
$\SpinC$ structure $\spincX$ which bounds $\Umet_Y$ and $\spinc$
respectively, we can find a $\SpinC$ connection ${\widetilde A}$ which
$dt + A$ in a collar neighborhood of its boundary. According
to~\cite{APSI}, the Dirac operator coupled to ${\widetilde A}$ acting
on $L^2(X)$ (given a cylindrical end $Y\times [0,\infty)$, and endowed
with the natural extension of ${\widetilde A}$) is a Fredholm
operator, and indeed the quantity 
\begin{equation}
\label{eq:DefCorrY}
\Corr_{\metY}(\spinc,A)=
\ind_{X}(\Dirac_{\widetilde A})-
\left(\frac{c_1({\widetilde A})^2-\sgn(X)}{8}\right),
\end{equation}
where $\ind_X(\Dirac)$ denotes the complex index of the Dirac operator,
and $c_1({\widetilde A})$, of course, denotes the Chern-Weil
representative of the first Chern class. 

According to the above proposition, if we take the difference between
$\CorrU(\spinc|U_0,\Umet_0,a_0)+\theta_{\met_t}(\spinc)+
\CorrU(\spinc|U_1,\Umet_1,a_1)$ and $\Corr_{\metY}(\spinc,A)$, we get a quantity
which is independent of the extending metrics $\Umet_0$, $\Umet_1$ and
the path $\met_t$.  In fact, we get something which is independent of
the Heegaard decomposition as well, and hence a topological invariant
of $Y$.

\begin{theorem}
\label{thm:WellDefd}
Let $\met_t$ be a $(U_0,U_1)$-allowable path, and the $\metY$ be the
metric formed from $\Umet_0$, $\met_t$, and $\Umet_1$, where
$\Umet_i\in\Met(U_i)$ are metrics which bound $\met_0$ and $\met_1$
respectively.  Then, for generic $a_0$ and $a_1$, the quantity
$$\NInv(\spinc)=\CorrU(\spinc|U_0,\Umet_0,a_0)+
\theta_{\met_t}(\spinc)+
\CorrU(\spinc|U_1,\Umet_1,a_1)
-\Corr_{\metY}(\spinc,A)$$
(where $A=A(\scale,T_0,T_1)$ and $\metY=\metY(\mu,T_0,T_1)$ 
for sufficiently large $\scale$, $T_0$,
and $T_1$)
is a topological invariant of $Y$.
\end{theorem}

The independence of $\NInv$ of the Heegaard decomposition
relies on the corresponding result for $\theta$, which was
established in~\cite{Theta}.  To state the result, use a connected sum
for paths of metrics. For a fixed metric on the torus $S^1\times S^1$,
a metric $\met$ over $\Sigma$ (which is flat in a neighborhood of the
connected sum point $p\in\Sigma$), and a real number $T>0$, let
$\met(T)$ denote the metric on $\Sigma'=\Sigma\#(S^1\times S^1)$
obtained by a connected sum with neck-length $T$. Similarly, for a
one-parameter family $\met_t$ over $\Sigma$, we let $\met_t(T)$ denote
the one-parameter family of metrics on $\Sigma'$ obtained in this
manner. Then, we have the following:

\begin{prop}
\label{prop:StabInvarTheta}
Let $U_0\#_{\Sigma} U_1$ be a Heegaard decomposition of $Y$, and let 
$U_0'\#_{\Sigma'}U_1'$ be the ``stabilized'' Heegaard decomposition; i.e.
$U_0'=U_0\#(S^1\times \CDisk)$, 
$\Sigma'=\Sigma\#(S^1\times S^1)$, 
$U_1'=U_1\#(\CDisk\times S^1)$. 
For a path of metrics on $\met_t$ and a neck-length $T$, let 
$\met_t(T)$ denote the path of metrics on $\Sigma'$ obtained by
forming the connected sum of metrics.
Given a generic $(U_0,U_1)$-allowable
path $\met_t$, for all sufficiently large $T$, $\met_t(T)$ is a
generic $(U_0',U_1')$-allowable path, and the moduli spaces are
diffeomorphic. Thus, 
$$\theta_{\met(T)}(\spinc)=\theta_{\met_t(T)}(\spinc).$$
\end{prop}

\begin{proof}
The allowability of $\met_t(T)$ and diffeomorphism statement for 
$\theta$ were
proved in Proposition~\ref{Theta:prop:Stabilization} of~\cite{Theta};
note that the hypothesis that $b_1(Y)>0$ was not used in the proof of
this fact. One can arrange for $\met_t$ and $\met_t(T)$ to be generic
simultaneously by varying the family $\met_t$ in a region $U\subset
\Sigma$ which does not contain the connected sum point. This 
statement is proved in
Proposition~\ref{prop:GenericMetrics}.
\end{proof}

The proof of Theorem~\ref{thm:WellDefd} is not difficult, given the
excisive properties of spectral flow, and the above
proposition. However, spelling out the precise form of excision needed
involves some notation which we give in the
Section~\ref{sec:Handlebodies}, so we defer the proof to the end of that
section. 

\section{Chambered metric invariants over Handlebodies}
\label{sec:Handlebodies}

Let $Y$ be a closed, oriented three-manifold equipped with a $\SpinC$
structure $\spinc$ with spinor bundle $W$ whose first Chern class is
torsion. Consider the space $\Conf(Y)$ of pairs $(\metY,A)$, where
$\metY$ is a metric over $Y$, and $A$ is a connection over $W$ (modulo
gauge). This space is homotopy equivalent to the torus
$\Torus{b_1(Y)}$. A function $$f\colon
\Conf(Y)-\{(\Umet,A)\big| \Ker
\Dirac_{(\Umet,A)}\neq 0\}
\longrightarrow \R$$ is said to be {\em a chambered metric
invariant} if
$$f(\metY',A')-f(\metY,A)=\SF(\Dirac_{(\metY,A)},\Dirac_{(\metY'A')}),$$
where $\SF((\metY,A),(\metY',A'))$ denotes the spectral flow of the
$\SpinC$ Dirac operator along any path in $\Conf(Y)$ which connects
the pairs $(\metY,A)$ and $(\metY',A')$; i.e. it is the intersection
number of the spectrum of the $\SpinC$ Dirac operator with the zero
eigenvalue.  Note that, since we have assumed that $c_1(W)$ is
torsion, the Atiyah-Singer index theorem guarantees that the spectral
flow is independent of the path.

Note that a chambered invariant is uniquely determined by its value on
any one pair $(\metY,A)$ over $Y$.  The quintessential chambered
invariant is the invariant $\Corr_{\metY}$ defined by
Equation~\eqref{eq:DefCorrY}. Our goal in this section is to construct
a chambered metric invariant for handlebodies. When working with
manifolds-with-boundary, special care must be taken to ensure that the
spectral flow used in the above definition makes sense.

Let $U$ be a handlebody which bounds $\Sigma$, and fix an
identification $\partial U\cong \Sigma$. Let $\met$ be a metric on
$\Sigma$. A metric $\Umet$ is said to {\em bound \met} if there is
a neighborhood of $\partial U$ which is isometric to
$(-1,0]\times\Sigma$ given the product metric $$(dt)^2+\met.$$ 
Note that a metric $\Umet$ which bounds a metric on $\Sigma$ can be
naturally extended to a cylindrical-end metric on 
the handlebody 
$$ U^+=U\cup_{\Sigma} \left([0,\infty)\times
\Sigma\right). $$ A metric $\Umet$ is said to be {\em product-like
near its boundary} if it bounds some metric $\met$ on its boundary. 
The relevant properties of the Dirac operator coupled to such a metric
are summarized in the following:

\begin{lemma}
\label{lemma:UDirac}
Let $\Umet$ be a metric on $U$ which bounds a $U$-allowable metric.
Then for all connections $A$ on the spinor bundle of the form $A=A_0 +
a$, where $A_0$ is reducible and $a$ is compactly supported, the
associated Dirac operator is Fredholm. Moreover, for generic such $A$,
the associated Dirac operator has no kernel.
\end{lemma}

\begin{proof}
The connection $A$ naturally induces a connection on the spinor bundle
of the boundary. More precisely, under the splitting 
$$W|\partial U = \SpinBunTwoP\oplus\SpinBunTwoM$$
into the $\pm 1$-eigenspaces of Clifford multiplication by $i$ times
the volume form of $\Sigma=\partial U$, the connection $A$ naturally
induces a connection $B$ on $\SpinBunTwoP$. If $\Tr F_A$ is compactly
supported in $U$, then the induced connection $B$ has normalized
curvature form.

The Dirac operator in a neighborhood of $\partial U$ takes the form
$$\DDt + \sqrt{2}\left(\begin{array}{cc}
0 & \DBar_B \\
\DBar_B^* & 0
\end{array}
\right) \colon \left(\begin{array}{c}
\SpinBunTwoP \\
\SpinBunTwoM
\end{array}
\right) \longrightarrow \left(\begin{array}{c}
\SpinBunTwoP \\
\SpinBunTwoM
\end{array}
\right). $$
According to Proposition~1.1 of~\cite{APSI}, this operator is Fredholm
if the kernels of $\DBar_B$ and $\DBar_B^*$ are trivial. Now, if $A$
is reducible, or even if it differs from a reducible by a compactly
supported one-form, then $B\in\Lag(U)$. Thus, the Fredholm condition
is guaranteed if $\Umet$ is $U$-allowable.

The genericity statement is an application of the Sard-Smale
theorem (see~\cite{SardSmale}). 
Let 
$$\ConfU(U)=\left\{(\Umet, A+a)\Bigg| 
\begin{array}{l}
\text{$\Umet$ bounds a $U$-allowable metric} \\
\Tr F_A = 0 \\
a\in \Omega^1_c(i\R)
\end{array}\right\},$$
and
$${\mathfrak M}=\{(\Umet,A,\Psi)\big|
(\Umet,A+a)\in\ConfU(U), 
\Dirac_{\Umet,A}\Psi=0, \|\Psi\|=1\}/\Maps(U,S^1).$$
With in Sobolev completions, ${\mathfrak M}$ is a Banach manifold
which is transversally cut out from $$\ConfU(U)\times
\{\Psi\big|\|\Psi\|=1\}/\Maps(U,S^1)$$ by the Dirac equation: i.e. if
$\phi$ were in its cokernel at $(\Umet,A,\Psi)$, then by varying the
spinor component, we would see that $\phi$ is
$\Dirac_{(\Umet,A)}$-harmonic. By varying the connection (indeed, in
any open set), we see that $\phi$ must vanish identically (by the
unique continuation principle).  It is easy to see that the projection
map from ${\mathfrak M}$ to $\ConfU(U)/\Maps(U,S^1)$ which forgets the
spinor is Fredholm of index $-1$.  It follows then from the Sard-Smale
theorem that for generic $(\Umet,A)\in \ConfU(U)$, there are no
harmonic spinors.
\end{proof}
 
The above lemma gives a space $\ConfU(U)$ of pairs $(\Umet,A)$ for
which the associated Dirac operator is a self-adjoint, Fredholm
operator. Thus, as in~\cite{APSIII}, the spectral flow along any path
in $\ConfU(U)$ is well-defined. Indeed, since the $\SpinC$ structure
on a handlebody comes from a spin structure, and $\ConfU(U)$ retracts
back to the space of flat connections modulo gauge, the spectral flow
of the Dirac operator between any two pairs in $\ConfU(U)$ is
well-defined and independent of the path joining them. (It is proved
in Lemma~\ref{Theta:lemma:SFHandlebody} of~\cite{Theta} that the
spectral flow around any loop in the space of flat connections over
the handlebody is trivial.)  Thus, a function 
$$f\colon \ConfU(U)-\{(\Umet,A)\big| \Ker
\Dirac_{(\Umet,A)}\neq 0\}
\longrightarrow \R$$
is said to be a {\em chambered invariant} if 
$$f(\Umet',A')-f(\Umet,A)=\SF(\Dirac_{\Umet,A},\Dirac_{\Umet',A'}).$$

The goal of this section is to define one canonical chambered
invariant $\CorrU$ for handlebodies of arbitrary genus, as we describe
shortly. We then spell out the properties of this invariant.  These
results rely on a splitting theorem for spectral flow for manifolds
with boundary, which in turn holds because we have a strong
non-degeneracy condition along corners and boundaries: in the form
which we require, the splitting principle can be seen to be a
straightforward consequence of the Fredholm theory developed by
M\"uller~\cite{MullerCorners}. After proving the various properties of
$\CorrU$, we use them to establish the stabilization invariance of the
normalized invariant $\NInv$, see Theorem~\ref{thm:WellDefd}.

\begin{figure}
\mbox{\vbox{\epsfbox{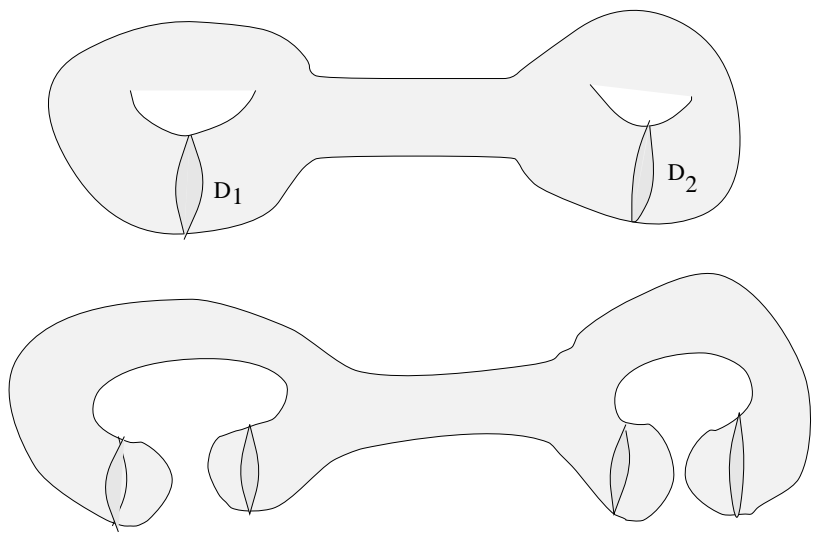}}}
\caption{\label{fig:AttachDisks}
Surgeries on a genus $2$ handlebody}
\end{figure}

To state the defining properties for $\CorrU$, we must describe how to
(metrically) perform surgeries on a handlebody.  
(Figure~\ref{fig:AttachDisks} gives a schematic 
illustration of this operation in the case where the handlebody has
genus two.)
Let $U$ be a
handlebody. Fix a collection of disjoint, embedded disks
$\{\CDisk_1,...,\CDisk_\NumDisks\}\subset U$ (whose boundaries are
embedded in the boundary of $\partial U$).  A metric $\Umet$ on $U$
which is product-like in a neighborhood of its boundary is said to be
{\em product-like in a neighborhood of the disks} if there is an
isometry from the metric product $$\coprod_{i=1}^\NumDisks
\CDisk_i \times [-1,1] \subset U,$$ which identifies the central
slice $\CDisk_i\times \{0\}$ with the $i^{th}$ disk, and where the
disks $\CDisk_i$ are endowed with a non-negative scalar curvature
metric which is product-like near its boundary. For positive real
numbers $T_1,...,T_\NumDisks$, let $\Umet(T_1,...,T_\NumDisks)$ denote
the metric stretched out normal to the disks: it is the metric
obtained by replacing the cylinder $\CDisk_i\times [-1,1]$ with the
elongated cylinder $\CDisk_i\times [-T_i,T_i]$. Moreover, we say that
the connection $A$ is product-like in a neighborhood of the disks
$\{\CDisk_1,...,\CDisk_\NumDisks\}$ if the trace of its curvature
vanishes in the product neighborhood and it is supported away from
$\partial U$. In this case, the connection has a canonical extension
$A(T_1,...,T_\NumDisks)$ over the stretched out handlebody.


Given a genus $g$ handlebody $U^g$, one can find $g$ embedded disks
$\{\CDisk_1,...,\CDisk_g\}$, whose complement in $U^g$ is homeomorphic
to a three-dimensional ball. 
Such a collection
$\{\CDisk_1,...,\CDisk_g\}$ is called a {\em complete set of attaching
disks} for $U^g$.  Note that if a metric $\Umet$ bounds a metric
$\met$ on its boundary, and is product-like in a neighborhood of a
complete set of attaching disks, then for all sufficiently large
$T_1,...,T_g$, the metric $\Umet(T_1,...,T_g)$ is $U^g$-allowable in a
neighborhood of its boundary. (Indeed, much of this section is modeled
on the corresponding results for $U$-allowable metrics in Section~2
of~\cite{Theta}.)

The complement $V$ of the attaching disks in $U^g$ is, strictly speaking,
a manifold-with-corners. We smooth out the corners to get a smooth
genus zero handlebody as follows.  Let $\HDisk$ be the
three-dimensional manifold-with-corners which is diffeomorphic to half
of the three-dimensional ball (i.e.  $\HDisk\cong
\{(x,y,z)\in\R^3\big| x^2+y^2+z^2\leq 1, z\geq 0$). This manifold has
two boundary components $\partial_0\HDisk$, $\partial_1\HDisk$, both
of which are two-dimensional disks, which meet normally along their
equator. The manifold obtained by attaching $2g$ copies of $\HDisk$ to
$V$ along its $2g$ faces is a smooth, genus zero handlebody.

If $\Umet$ is a metric on $U^g$ which is product-like in a
neighborhood of the $g$ attaching disks, the genus zero handlebody
inherits a metric which depends on $g$ parameters, denoted
$\Umet_0(T_1,...,T_g)$. This metric depends also on a choice 
of a metric $\PCMet$ on  $\HDisk$ of non-negative sectional
curvature, which in a neighborhood of both of its boundaries is
isometric to $\CDisk\times [0,\epsilon)$. We construct such a 
metric in Lemma~\ref{lemma:NNCPiece}. Given this result, we let
$\Umet_0(T_1,...,T_g)$ denote the metric obtained by removing
$(-1,1)\times \CDisk_i$ (for $i=1,...,g$) from $U^g$, attaching
solid cylinders
$[0,T_i]\times \CDisk$ along the $2g$ new boundary components, and
then capping off with $2g$ copies of $\HDisk$, endowed with the metric
$\PCMet$. Formally, $\Umet_0(T_0,...,T_g)$ is the metric inherited
from the description:
$$\Big(U^g - \bigcup_{i=1}^g [-1,1]\times \CDisk_i \Big)
\cup_{\{\pm 1\}\times
\CDisk_i=\{0\}\times \CDisk} 
\Big(\left[0,T_i\right]\times \CDisk\Big)
\cup_{\{T_i\}\times
\CDisk = \partial_1 B_{\pm i}} 
\Big(B_{\pm i}\Big),$$
where $B_{\pm i}$ are $2g$ copies of $\HDisk$.
Note that $\Umet_0(T_1,...,T_g)$ is a metric on the 
genus zero handlebody (i.e. a three-ball), which is product-like
near its boundary. 

Moreover, if we have a $\SpinC$ connection $A$ on $U^g$ which is
product-like in a neighborhood of the disks, then it has a natural
extension to the surgered manifold, obtained by extending it over the
three-balls $\HDisk$ to have traceless curvature. We denote the
resulting connection by $A_0(T_1,...,T_g)$. Before stating the
definition of $\CorrU$, we pause to construct the metric $\PCMet$ on
the three-dimensional half-ball used in the above construction.

\begin{lemma}
\label{lemma:NNCPiece}
There is a metric $\PCMet$ on $\HDisk$ with everywhere non-negative sectional
curvatures, which is product-like in a neighborhood of its boundaries,
and whose boundary is a union of two copies of $\CDisk$, meeting at a
corner.
\end{lemma}

\begin{proof} 
Fix some constant $0\leq \epsilon < \frac{1}{4}$, and fix a smooth,
non-decreasing function $$\psi\colon [0,1]\longrightarrow [0,1]$$ with
$$
\psi(t)=\left\{\begin{array}{ll}
\frac{3\epsilon}{2} &{\text{if $t\leq \epsilon$}} \\
t &  {\text{if $t\geq 2\epsilon$}} 
\end{array}\right.
$$
Then, the hypersurface $$S^3 \cong \{(x_1,x_2,x_3,x_4)\in \R^4\big|
\sum_{i=1}^4 \psi(x_i)^2 = 1\}$$ inherits a metric from $\R^4$. The
region where $x_3\geq 0$ and $x_4 \geq 0$ is diffeomorphic to half of
the three-ball (its two boundaries, are the loci where $x_3$ and $x_4$
vanish respectively). Since the function $\psi$ is constant for small
values, the metric is easily seen to respect the corners.  Moreover,
the metric is easily seen to have all non-negative sectional
curvatures (see for example~\cite{DoCarmo}), as the hypersurface is
locally described as a graph of $\sqrt{1-\psi(x_1)^2 -
\psi(x_2)^2-\psi(x_3)^2}$ (after possibly renumbering the four
variables), which is clearly a convex function.
\end{proof}

By doubling the metric $k_0^+$ above, we obtain a cylindrical-end
metric on the three-ball $D^3$ with non-negative scalar curvature
(note that this can be connected to any ``standard'' cylindrical-end
metric on the three-disk through metrics of non-negative scalar
curvature). 

\begin{prop}
\label{prop:CorrU}
There is a unique chambered invariant $\CorrU$ with the following 
properties:
\begin{enumerate}
\item If $g=0$, and we endow $U$ with the metric obtained by doubling $k_0^+$
(and a connection with traceless curvature), then $\CorrU(U)=0$.
\item If $g=1$ and we endow $U$ with a product 
metric of the form $S^1\times \CDisk$ (and a connection with traceless
curvature), then $\CorrU(U)=0$.
\item If $(\Umet,A)$ is product-like in a neighborhood of $g$ attaching disks, then 
$$\CorrU(\Umet(T),A(T))=\CorrU(\Umet(T'),A(T')),$$
for all sufficiently large $T,T'$.
\item If $(\Umet,A)$ is product-like in a neighborhood of $g$ attaching disks, then 
\begin{equation}
\label{eq:DefRel}
\CorrU(\Umet(T),A(T))=\CorrU(\Umet_0(T),A_0(T))
\end{equation}
for all sufficiently large $T$.
\end{enumerate}
\end{prop}

This invariant is additive for boundary connected sum, in the following sense. Let
$U^{g_1}$ 
and $U^{g_2}$ be a pair of handlebodies, and fix points
$p_1$ and $p_2$ in $\partial U^{g_1}$ and $\partial U^{g_2}$
respectively. Fix metrics $\Umet_1$ and $\Umet_2$ for which a neighborhood of 
the $p_1$ and $p_2$, are isometric to the standard piece
$\HDisk - \partial_1\HDisk$. Then, we can form the ``boundary connected sum''
$$U^{g_1}\#_T U^{g_2} = \left(U^{g_1}-\HDisk\right) \cup \left([-T,T]\times
\CDisk\right) \cup \left(U^{g_2}-\HDisk\right).$$
This is a genus $g_1+g_2$ handlebody, endowed with a metric, denoted
$\Umet_1\#_T\Umet_2$, which bounds the connected sum metric
$\left(\partial U^{g_1}\right)\#_T \left(\partial U^{g_1}\right)$.
Moreover, if $A_1$ and $A_2$ are connections whose curvature is
traceless over the standard pieces, then the connections can be
naturally extended, as well.

\begin{prop}
\label{prop:Additivity}
The invariant $\CorrU$ is additive under boundary connected sum, in
the sense that for pairs $(\Umet_1,A_1)$ and $(\Umet_2,A_2)$ on
$U^{g_1}$ and $U^{g_2}$, there is a $T_0>0$, so that for all $T\geq
T_0$, $$\CorrU((\Umet_1,A_1)\#_{T}(\Umet_2,A_2))=
\CorrU((\Umet_1,A_1)\#_T (\PCMet,\PCA))
+\CorrU((\Umet_2,A_2)\#_{T}(\PCMet,\PCA)),$$ where $(\PCMet,\PCA)$ is
a pair over $\HDisk$ where $\PCMet$ is as in
Lemma~\ref{lemma:NNCPiece}, and $\PCA$ has traceless curvature.
\end{prop}

This invariant is compatible with the correction term $\Corr$ for
$S^3$ with its genus zero or one Heegaard decompositions, in the
following sense:

\begin{prop}
\label{prop:CompareCorrections}
Let $S^3=U_0\cup_{\Sigma} U_1$ be a standard genus zero or one
Heegaard decomposition of $S^3$.  Let $\met_t$ be a path of metrics on
the Heegaard surface (sphere or two-torus) $\Sigma$, and let
$\Umet_0$, $\Umet_1$ be metrics over the genus zero or one
handlebodies $U_0$ and $U_1$ which extend $\met_0$ and $\met_1$
respectively. Then,
$$\Corr\big((\Umet_0)\#(\met_t)\#(\Umet_1)\big)=\CorrU(\Umet_0)+\CorrU(\Umet_1).$$
\end{prop}

The construction of $\CorrU$, and the proof of its various properties,
rests a splitting (or excision) principle for spectral flow, which we
will presently outline.  This excision principle involves degenerating
the handlebodies normal to embedded disks, and the objects one
encounters under such degenerations are
manifolds-with-corners. 
Formally, we have the following:

\begin{defn}
A {\em truncated handlebody} $V$ is a smooth three-manifold with
corners which is homeomorphic to a handlebody, and whose codimension
one boundary consists of a surface-with-boundary $\Surf$, the {\em
bounding surface}, and a collection of disjoint disks
$\{\CDisk_1,...,\CDisk_\NumDisks\}$, the {\em faces}. The bounding
surface $\Surf$ meets the faces normally along the boundary.
\end{defn}

Examples of truncated handlebodies include the product $\CDisk\times
[0,1]$, the half-ball $\HDisk$, and the complement of a collection of
embedded disks in the genus $g$ handlebody.
%
%

It is useful to describe the product structure near the boundaries of
a truncated handlebody $V$ in detail. To cover the faces, we
have a diffeomorphism $$\Phi\colon
\bigcup_{i=1}^\NumDisks
\CDisk_i\times (-1,0] \subset V$$ 
onto a neighborhood of this boundary region for $V$. In turn, 
a neighborhood of the boundary of the union of disks admits an 
identification
$$\phi\colon \bigcup_{i=1}^{\NumDisks}(-1,0]\times S^1_{i}
\longrightarrow \bigcup_{i=1}^{\NumDisks}\CDisk_{i}.$$ 
The remaining boundary
region for $V$ is a surface $F$ of genus zero with $\NumDisks$ boundary
circles, and we have a diffeomorphism
$$\Psi\colon
F\times (-1,0]\subset V.$$
onto this boundary region.
A neighborhood fo the boundary of $F$ admits an identification
$$\psi \colon \bigcup_{i=1}^{\NumDisks}S^1_{i}\times
(-1,0]
\longrightarrow F.$$ We can require that all these identifications be
compatible, in the sense that the following maps commute:
$$\begin{CD}
\bigcup_{i=1}^{\NumDisks}(-1,0]\times S^1_{i}\times (-1,0]
@>{\phi\times \Id}>>
\bigcup_{i=1}^{\NumDisks}\CDisk_i\times (-1,0] \\
@V{\Id\times \psi}VV  @V{\Phi}VV \\
(-1,0]\times F @>{\Psi}>> V
\end{CD}.
$$

A metric $\Vmet$ is said to be {\em product-like in a neighborhood of
its boundaries} if the above identifications $\Phi, \Psi, \phi, \psi$
are all isometries, where the domains of the maps are all given
product metrics (and their ranges are given metrics induced from
$\Vmet$). 

A truncated handlebody $V$ can be naturally completed to a
three-manifold $V^+$ without boundary by attaching $[0,\infty)\times
\Surf$ along the bounding surface, solid cylinders $\CDisk_i\times
[0,\infty)$ along the faces, and a region $[0,\infty)\times
S^1_i\times [0,\infty)$ along the corners $\{S^1_i\}$.  If a metric
$\Vmet$ is product-like in a neighborhood of its
boundaries, it can be naturally extended to a complete metric, by
giving all these standard pieces product metrics in a compatible
manner. In particular, this compatibility ensures that in the region
$$\left([0,\infty)\times \Surf \right)\cup \left([0,\infty)\times
S^1_i\times [0,\infty)\right),$$ the metric $\Vmet^+$ is isometric to a
product of $[0,\infty)$ with the cylindrical completion of the
bounding surface $$\Surf^+=\Surf_{\partial\Surf=\bigcup
S^1_{i}}\bigcup S^1_i\times [0,\infty)$$

As always, we will consider the Dirac operator with respect to such a
metric, coupled to a $\SpinC$ connection $A$ with traceless
curvature. Note that the connection $A$ can be naturally extended to a
connection $A^+$ on $V^+$ in such a manner that the curvature remains
traceless.

\begin{defn}
The pair $(\Vmet,A)$ is said to be {\em strongly non-degenerate on the
boundary} if the restriction $(\Vmet, A)$ to the bounding surface
$\Surf$ has trivial kernel with APS boundary conditions or,
equivalently, if the induced Dirac operator induced on any of the
attached slices $\Surf^+\subset V^+$ has trivial $L^2$ kernel.
\end{defn}

Note that there is another product region ``at infinity'', the region 
$$\CDisk^+\times [0,\infty) = \left(\CDisk_i\times [0,\infty)\right)\cup
\left([0,\infty)\times S^1_i\times [0,\infty)\right).$$
The induced Dirac operator on these attached slices $\CDisk^+$
automatically has trivial $L^2$ kernel, since its kernel is naturally
identified with the harmonic spinors on the two-sphere (see
Proposition~\ref{Theta:prop:LTwoCohom} of \cite{Theta}).

In~\cite{MullerCorners}, M\"uller considers Dirac operators on
manifolds-with-corners, which satisfy a certain non-degeneracy
hypothesis along its corners (see also~\cite{MazzeoMelrose}). Specializing some of his results to the
case of truncated handlebodies, we obtain a Fredholm and exponential
decay result. To spell out the exponential decay result, let
$V_T\subset V^+$ denote the subset obtained by attaching subsets
$\Surf\times [0,T]$, $[0,T]\times S^1_i \times [0,T]$, and
$\CDisk_i\times [0,T]$ to $V$.

\begin{prop}
\label{prop:Mueller}
Let $V$ be a truncated handlebody.  Suppose that $(\Vmet,A)$ is
strongly non-degenerate on the boundary. Then, the Dirac operator
coupled to $A^+$ induces a Fredholm operator on $L^2(V^+)$.  In
particular, there is a real number $\epsilon>0$ with the property that
the spectrum of the Dirac operator in the range $(-\epsilon,\epsilon)$
is discrete. Moreover, eigenvectors in this range enjoy an exponential
decay property: there are constants $C,c>0$ with the property that for
each $\lambda$-eigenvector $\Phi$ of $\Dirac_{A^+}$ for $\lambda\in
(-\epsilon,\epsilon)$, $$\int_{V^+-V_T} |\Phi|^2 \leq C
e^{-cT}\int_{V^+}|\Phi|^2.$$
\end{prop}

\begin{proof}
Both statements are proved in~\cite{MullerCorners}: the Fredholm
paramatrix is constructed in the proof of Proposition~2.8, and the
exponential decay estimate is Proposion~2.19 of that reference.
\end{proof}

Given the Fredholm paramatrix and exponential decay, the usual
splicing techniques give a splitting principle for spectral
flow. Suppose that $V$, $V'$ be a pair of truncated handlebodies, and
let $\{\CDisk_1,...,\CDisk_\NumDisks\}$,
$\{\CDisk_1',...,\CDisk_\NumDisks'\}$ be a collection of (not
necessarily all) faces in $V$
and $V'$ respectively. Then, we can form a new truncated handlebody
$$V\#_T V' = \left(V\right)\cup_{\CDisk_i = \CDisk_i\times \{-T\}}
\left(\bigcup_{i=1}^\NumDisks \CDisk_i\times [-T,T]\right)
\cup_{\CDisk_i\times \{T\}=\CDisk_i'} \left(V'\right).$$

If $(\Vmet,A)$ and $(\Vmet',A')$ is are strongly non-degenerate pairs
on $V$ and $V'$, discrete spectrum and exponential decay
considerations on the boundaries show that for all $T$ sufficiently
large, the induced pair $(\Vmet,A)\#_T(\Vmet',A')$ is strongly
non-degenerate on $V\#_T V'$. Indeed, we have the following
direct consequence of Proposition~\ref{prop:Mueller}:

\begin{prop}
\label{prop:SplittingPrinciple}
Let $(\Vmet_t, A_t)$ and $(\Vmet_t',A_t')$ be a pair of one-parameter
families of strongly non-degenerate data on $V$ and $V'$, then for all
sufficiently long tube-lengths, the spectral flow of $V\#V'$ is the
sum of the spectral flows of the pieces. More precisely, if the kernels
of the Dirac operators of $\Dirac_{(\Vmet_0,A_0)}$,
$\Dirac_{(\Vmet_1,A_1)}$,  $\Dirac_{(\Vmet_0',A_0')}$,
$\Dirac_{(\Vmet_1',A_1')}$ are all trivial, then there is a real $T_0$
so that for all $T\geq T_0$, the kernels of
$\Dirac_{(\Vmet_0,A_0)\#_T(\Vmet_0',A_0')}$ and 
$\Dirac_{(\Vmet_1,A_1)\#_T(\Vmet_1',A_1')}$ are trivial, and indeed
$$\SF((\Vmet_t,A_t)\#_T(\Vmet_t',A_t'))=\SF(\Vmet_t,A_t)+\SF(\Vmet_t',A_t').$$
\end{prop}

Having set up the splitting principle for spectral flow, we turn our
attention to the existence and uniqueness for $\CorrU$:

\vskip.2cm
\noindent{\bf{Proof of Proposition~\ref{prop:CorrU}}}.
Suppose $g=0$. Let $\Umet_0$ be a metric obtained by doubling $\PCMet$. 
Let $A_0$ be the connection on the
corresponding spinor bundle with traceless curvature. Given any pair
$\data=(\Umet,A)\in \ConfU(U)$, define
$$\CorrU(\data)=\SF(\PCdata,\data),$$
where $\PCdata$ consists of the metric $\Umet_0$ and the reducible 
connection $A_0$.
This is by definition a
chambered invariant.  Hence we have existence
and uniqueness when $g=0$.

Suppose now that $g>0$, and let $\data=(\Umet,A)\in\ConfU(U)$ be a
pair which is product-like in a neighborhood of $g$ attaching disks
$\{\CDisk_1,...,\CDisk_g\}$. First, we show that for all sufficiently
large $T$ the right hand side of Equation~\eqref{eq:DefRel}
stabilizes. Let $V$ denote the truncated handlebody obtained by
removing the product neighborhoods of the attaching disks from $U^g$.
For generic $\data$, the kernel of the Dirac operator on $V$ has no
kernel (this follows from the Fredholm property from 
Proposition~\ref{prop:Mueller}, together with the 
proof of Lemma~\ref{lemma:UDirac}), 
so, by the splitting principle for the zero-modes, there is a
$T_0$ so that for all $T\geq T_0$, the right hand side of
Equation~\eqref{eq:DefRel} for $\CorrU(\data)$ does not depend on
$T$. So we let Equation~\eqref{eq:DefRel} be the definition of
$\CorrU(k(T),A(T))$. Using spectral flow, this specifies $\CorrU$ for
any pair $(k,A)\in\ConfU(U)$. We have to show that $\CorrU$ is independent
of the choice of product-like metric we started with, and, indeed, the
choice of attaching disks. 

To this end, let  $\data=(\Umet,A)$ and $\data'=(\Umet',A')$ be
product-like in a neighborhood of the same $g$ attaching
disks. Then by the splitting principle for spectral flow, for all
sufficiently large $T$
$$\SF(\data(T),\data'(T))=\SF(\data_0(T),\data'_0(T)).$$
It follows that $\CorrU$ is independent of the choice of initial $\data$.

\begin{figure}
\mbox{\vbox{\epsfbox{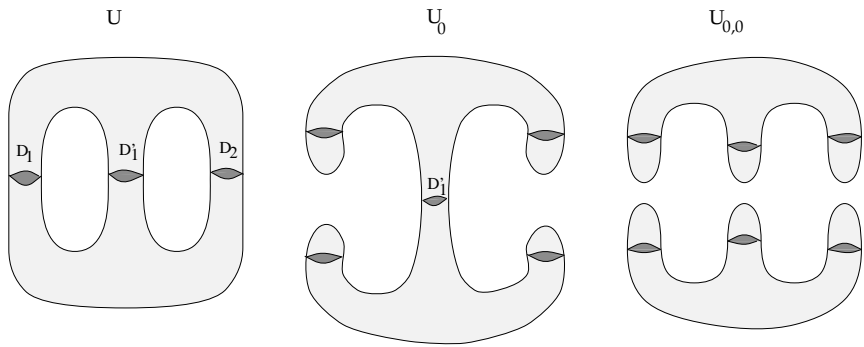}}}
\caption{
    \label{fig:DoubleSurger}
A genus $2$ handlebody $U$ with two complete sets of 
attaching disks $\{D_1,D_2\}$ and $\{D_1',D_2\}$. The
manifold $U_0$ obtained by surgering along
the first set of attaching disks is a three-ball, and the manifold
$U_{0,0}$ obtained by surgering along all the attaching disks
is a  disjoint union of
two three-balls.}
\end{figure}

Next, we verify that the definition of $\CorrU$ is independent of the
particular choice of attaching disks when $g>0$. Suppose that
$\data\in\ConfU(U)$ is product-like normal to $g+1$ disks
$\CDisk_1',\CDisk_1,\CDisk_2,...,\CDisk_g$, so that
$\{\CDisk_1,...,\CDisk_g\}$ and $\{\CDisk_1',\CDisk_2,...\CDisk_g\}$
are both complete sets of attaching disks. We would like to show that
the value of $\CorrU$ obtained by surgering out the first set is the
same as that obtained by surgering out the second set. To see this, we
prove that the invariant agrees with a quantity obtained by surgering
out all $g+1$ disks simultaneously.  More precisely, let
$\data_0(T',T)$ denote the pair on the genus zero handlebody obtained
by inserting a cylinder $[-T',T]\times\CDisk$ about $\CDisk_1'$ and
surgering out the $\{\CDisk_1,...,\CDisk_g\}$, and let
$\data_{0,0}(T',T)$ denote the metric obtained by surgering all $g+1$
disks (using tubelength $T'$ around $\CDisk_1'$ and $T$ around all the
others). Note that the pair $\data_{0,0}(T',T)$ lives on a union of
two disjoint genus zero handlebodies. (See
Figure~\ref{fig:DoubleSurger} for an illustration.)


By the splitting principle, the Dirac operator is generically
invertible for all sufficiently long tubelengths (along all $g+1$
disks) in both $U$ and the doubly-surgered $U_{0,0}$.  Thus, the
quantity
$$Q=\SF(\data,\data(T',T))+\SF(\data_{0,0}(T',T),\PCdata\coprod\PCdata)$$
is independent of the particular $T',T$, provided that both are larger
than some constant $T_0$.  Our aim is to show that $Q$, in which
both sets of attaching disks play the same role, agrees with $\CorrU$
calculated using the disks $\{\CDisk_1,...,\CDisk_g\}$.  We compare
$Q$ with $$P=\SF(\data,\data(1,T))+\SF(\data_0(1,T),\PCdata),$$ the
quantity obtained by calculating $\CorrU$ using the disks
$\{\CDisk_1,...,\CDisk_g\}$.  Now
\begin{eqnarray*}
P-Q &=&
\SF(\data(T',T),\data(1,T)) + \SF(\data_0(1,T),\PCdata) +
\SF(\PCdata\coprod\PCdata, \data_{0,0}(T',T)).
\end{eqnarray*}
But $$\SF(\data_{0,0}(T',T),\PCdata\coprod\PCdata)
= \SF(\data_{0}(T',T),\PCdata),$$ since taking
$T$ to be sufficiently large, both terms can be identified with the
spectral flow from the complement of $\CDisk_1'$ in
$\data_{0}(1,T)$ to the disjoint union of two half-balls
$\HDisk\coprod \HDisk$ (endowed with the metric $k_0$ and connections
with traceless curvature). Substituting back, we see that $P=Q$. Since
we can switch the roles of $\CDisk_1$ and $\CDisk_1'$ (without
changing $Q$), we have shown that $\CorrU$ is independent of the
attaching disks.

Note that $\CorrU(S^1\times \CDisk)$ is independent of the length of
the $S^1$ factor, since the spectral flow between two such metrics
vanishes. So, by definition $\CorrU(S^1\times \CDisk)$ agrees with the
$\CorrU$ for the metric on $D^3$ obtained by surgering out the
attaching disk. This metric is precisely $\PCMet\cup [-T,T]\times
\CDisk \cup \PCMet$ -- so $\CorrU$ of it vanishes.
\qed
\vskip.2cm

\vskip.2cm
\noindent{\bf{Proof of Proposition~\ref{prop:Additivity}.}} The proof
follows in the same manner as the proof that $\CorrU$ is independent
of the attaching disks. Consider the attaching disks for $U^{g_1}$ and
$U^{g_2}$, and the disk used for the boundary connected sum. Then pull out all 
$g_1+g_2+1$ disks simultaneously.
\qed
\vskip.2cm

\vskip.2cm
\noindent{\bf{Proof of Proposition~\ref{prop:CompareCorrections}.}}
By the splitting theorem for spectral flow on the closed manifold
$S^3$, we can reduce the proposition to a model case.  Suppose that
$\Umet_0$, $\met_t$, and $\Umet_1$ are metrics on $U_0$, $\Sigma$, and
$U_1$ for which the proposition is known, and 
let $\Umet_0'$,  $\met_t'$, and $\Umet_1'$ denote arbitrary (compatible) 
metrics on
$U_0$, $\Sigma\times[0,1]$, and $U_1$. Then,
\begin{eqnarray*}
\Corr\left((\Umet_0')\#(\met_t')\#(\Umet_1')\right)
 -
\Corr\left((\Umet_0)\#(\met_t)\#(\Umet_1)\right)
&=&
\SF\left((\Umet_0)\#(\met_t)\#(\Umet_1),
(\Umet_0')\#(\met_t')\#(\Umet_1')\right) \\
&=& \SF(\Umet_0,\Umet_0') + \SF(\met_t,\met_t') +
\SF(\Umet_1,\Umet_1').
\end{eqnarray*}
Note now that for $i=0,1$,
$\SF\left(\Umet_i,\Umet_i'\right)=\CorrU(\Umet_i')-\CorrU(\Umet_i)$.
Moreover, it follows from Proposition~\ref{prop:SFCyl} that
$\SF(\met_t,\met_t')\equiv 0$. Thus, the proposition is established
once it is established for a model triple of metrics.

We consider the case where $g=1$. Let $\Umet_0$ and $\Umet_1$ be a
pair of metrics the of the form $\CDisk\times S^1$, where the disk is
endowed with a non-negative scalar curvature metric, and the $S^1$
factor has the same length as the boundary of $\CDisk$; let $\met_t$
denote the constant family of metrics. We know that
$\CorrU(\Umet_0)=\CorrU(\Umet_1)=0$, and need to show that for all
sufficiently large $T$, $\CorrY(\Umet_0\#_T\Umet_1)=0$.

We connect the standard, round metric on $S^3$ (for which it is easy
to see that $\CorrY=0$, since it bounds a metric on the four-ball with
non-negative sectional curvatures) with a metric of the form
$\Umet_0\#_T
\Umet_1$ through a path of metrics with non-negative scalar curvature,
to show that the correction terms $\CorrY$ agree. To this end, let
$$\psi_s\colon [0,1]\rightarrow [0,1]$$ be a smooth, one-parameter
family of smooth functions, depending smoothly on a parameter
$s\in[0,1]$, with following properties:
\begin{list}
{{(\arabic{bean})}}{\usecounter{bean}\setlength{\rightmargin}{\leftmargin}}
\item $\psi_s(t)=t$ for $t<1-s$,
\item $\DDt\psi_s(t)\equiv 0$ for $t>1-s+\epsilon$,
\item $\DDt\psi_s(t)\geq 0$ for all $t$,
\item $\frac{d^2\psi_s}{dt^2}(t)\leq 0$ for all $t$.
\end{list}
For example, if $f\colon \R\rightarrow [0,1]$ is a smooth,
non-increasing, non-negative function with $f(t)\equiv 1$ for
$t<0$, $f(t)\equiv 0$ for $t>\epsilon$, then we can let
$\psi_s(t)=\int_{0}^t f(x+s-1)dx$.

Let $(r,\theta,\phi)$ denote coordinates on $(0,1)\times S^1\times
S^1$. The standard three-sphere can be obtained from this space by
attaching two circles $0\times 0\times S^1$ and $1\times S^1\times 0$
``at infinity'' (i.e. at $r=0$ and $r=1$). 
Moreover, for all $s<1$, the metric on $(0,1)\times
S^1\times S^1$ given by $$g_s=\frac{dr^2}{1-r^2} + \psi_s(r)^2
d\theta^2 +
\psi_s(\sqrt{1-r^2})^2d\phi^2$$ extends over the two circles at $r=0$,
$r=1$ to give a smooth metric on $S^3$, since in a neighborhood of
those regions, $\psi_s(r)\equiv r$ and $\phi_s(\sqrt{1-r^2})\equiv
\sqrt{1-r^2}$. Note that when $s=0$, the above metric agrees with the standard round metric on $S^3$. 
By Cartan's method of moving frames (see~\cite{Spivak}), it is easy to
see that the connection matrix of the Levi-Civita connection is given
by: $$\left(\begin{array}{ccc} 0 &
\sqrt{1-r^2}\psi_s'(r)d\theta & -r \psi_s'(\sqrt{1-r^2}) d\phi \\
-\sqrt{1-r^2}\psi_s'(r)d\theta & 0 & 0 \\ r \psi_s'(\sqrt{1-r^2})
d\phi & 0 & 0
\end{array}
\right);$$
and hence the curvature matrix is:
$$\left(\begin{array}{ccc}
0 & -A(r)dr\wedge d\theta & -B(r)dr\wedge d\phi \\
A(r)dr\wedge d\theta & 0 & 0 \\
B(r)dr\wedge d\phi & 0 & 0 
\end{array}
\right),$$
where $A_s(r)=\frac{\psi_s'(r)r}{\sqrt{1-r^2}}-
\psi_s''(r)\sqrt{1-r^2}$ and $B_s(r)=\psi_s'(\sqrt{1-r^2})-\frac{r^2
\psi_s''(r)}{\sqrt{1-r^2}}$

Thus, the sectional curvatures are always non-negative. 

When $s>\OneHalf$, then $g_s$ extends to a decomposition of $S^3$ as a
union of $S^1\times \CDisk$ with $\CDisk\times S^1$ endowed with the
product metric. In particular, the metric is product-like in the
region where $r$ is in a neighborhood of $\OneHalf$.

Thus, the proposition follows.
\qed
\vskip.2cm

\vskip0.3cm
\noindent
{\bf{Proof of Theorem~\ref{thm:WellDefd}.}} 
In view of Proposition~\ref{prop:ExtendToY}, the invariant
$\NInv(\spinc)$ can depend only on the Heegaard decomposition, not on
the metrics used in its definition. Topological invariance then
amounts to showing that $\NInv(\spinc)$ remains unchanged under
stabilization. But this is a consequence of the
excision property of indices, together with the stabilization of 
$\Inv(\spinc)$ (see Proposition~\ref{prop:StabInvarTheta}).

Suppose $U_0\cup_{\Sigma} U_1$ is a Heegaard decomposition
and 
$$\left(U_0\#\left(S^1\times\CDisk\right)\right) \cup_{\Sigma\#\left(S^1\times S^1\right)}
\left(U_1\#\left(\CDisk\times S^1\right)\right)$$
is its stabilization, then the difference $\Delta$ between the corresponding
invariants
$\NInv$ (which we would like to show vanishes) is given by:
\begin{eqnarray}
\Delta&=&
\CorrU\left(U_0\#_T\left(S^1\times \CDisk\right)\right)
+\theta\left(\met_t\#_T\left(S^1\times S^1\right)\right)
+\CorrU\left(U_1\#_T\left(\CDisk\times S^1\right)\right) \nonumber \\
&&-\Corr\left(\left(U_0\#_T\left(S^1\times \CDisk\right)\right)\#
\left(\met_t\#_T\left(S^1\times
S^1\right)\right)\# \left(U_1\#_T\left(\CDisk\times
S^1\right)\right)\right) \nonumber \\
&&-\CorrU\left(U_0\#_T\TBall\right)
-\theta\left(\met_t\#_T\left(S^2\right)\right)
-\CorrU\left(U_1\#_T\TBall\right) \nonumber \\
&&+\Corr\left(\left(U_0\#_T\TBall\right)\#
\left(\met_t\#_T\left(S^2\right)\right)\#
\left(U_1\#_T\left(\TBall\right)\right)\right). \label{eq:DifferenceOne}
\end{eqnarray}
Note that on $\Sigma\# (S^1\times S^1)$ we are using an allowable path
induced from an allowable path on $\Sigma$, as in
Proposition~\ref{prop:StabInvarTheta}; also on the handlebodies, we
are using the metric arising from the boundary connected sum.
The difference in the terms using $\Corr$ is a
spectral flow (by its chambered nature); moreover by the excision
principle for spectral flow, we get:
\begin{eqnarray}
\lefteqn{\Corr\left(\left(U_0\#_T\left(S^1\times \CDisk\right)\right)\#
\left(\met_t\#_T\left(S^1\times
S^1\right)\right)\# \left(U_1\#_T\left(\CDisk\times
S^1\right)\right)\right)} \nonumber \\
\lefteqn{-\Corr\left(\left(U_0\#_T\TBall\right)\#
\left(\met_t\#_T\left(S^2\right)\right)\# \left(U_1\#_T\left(\TBall\right)\right)\right)}
\nonumber \\
&=&
{\Corr\left(\left(\TBall\#_T\left(S^1\times \CDisk\right)\right)\#
\left(S^2\#_T\left(S^1\times
S^1\right)\right)\# \left(\TBall\#_T\left(\CDisk\times
S^1\right)\right)\right)} \nonumber \\
&&-\Corr\left(\left(\TBall\#_T\TBall\right)\#
\left(S^2\#_T\left(S^2\right)\right)\#
\left(\TBall\#_T\left(\TBall\right)\right)\right) \nonumber \\
&=& {\Corr\left(\left(\TBall\#_T\left(S^1\times \CDisk\right)\right)\#
\left(S^2\#_T\left(S^1\times
S^1\right)\right)\# \left(\TBall\#_T\left(\CDisk\times
S^1\right)\right)\right)} 
 \label{eq:DifferenceTwo}
\end{eqnarray}
i.e. we have excised out $Y$ and replaced it by a three-sphere (and
used the positivity of the scalar curvature on the second term). Now,
the compatibility of $\CorrU$ and $\Corr$
(Proposition~\ref{prop:CompareCorrections}) allows us to conclude that
\begin{eqnarray}
\lefteqn
{\Corr\left(\left(\TBall\#_T\left(S^1\times \CDisk\right)\right)\#
\left(S^2\#_T\left(S^1\times
S^1\right)\right)\# \left(\TBall\#_T\left(\CDisk\times
S^1\right)\right)\right)}  \nonumber \\
&=& \CorrU(\TBall\#_T(S^1\times \CDisk)) +
\CorrU(\TBall\#_T(\CDisk\times S^1)) \label{eq:DifferenceThree}.
\end{eqnarray}
Combining Equations~\eqref{eq:DifferenceOne}, \eqref{eq:DifferenceTwo}, and
\eqref{eq:DifferenceThree} (and using the additivity of $\CorrU$, see Proposition~\ref{prop:Additivity}), we get that
$$\Delta=\theta(\met_t\#(S^1\times S^1))-\theta(\met_t\#(S^2)).$$
By the stabilization invariance of $\theta$ (see
Proposition~\ref{prop:StabInvarTheta}), this
implies that $\Delta=0$, as required.
\qed
\vspace{0.2in}


\section{Transversality of the Theta Divisor}
\label{sec:GenMet}

The aim of this section is to prove a ``generic metrics'' result for
the Abel-Jacobi map. Indeed, we show that for a fixed divisor
$D\in\Sym^{g-1}(\Sigma)$, the map from the space of metrics to the
Jacobian $$h\mapsto \Theta_{h}(D)$$ is a submersion onto the Jacobian,
provided that the genus $g$ of the Riemann surface is greater than
one.  Strictly speaking, to view this as a map into one, fixed space,
we fix a spin structure over $\Sigma$, and hence an identification
$\Jac_\met\cong H^1(\Sigma;S^1)$.

This technical result was used in our definition of $\Inv$: it
guarantees that the moduli space whose count defines $\Inv$ does
consist of isolated points. (It was also used in
Section~\ref{Theta:sec:WallCross} of~\cite{Theta}, in the derivation of
a ``wall-crossing'' formula.)

\begin{prop}
\label{prop:GenericMetrics}
Let $\Sigma$ be a surface of genus greater than one. Then, for each
divisor $D\in\Sym^{g-1}(\Sigma)$, the map
$$\Met(\Sigma)\longrightarrow H^1(\Sigma;S^1)$$ given by $\met\mapsto
\Theta_{\met}(D)$ is a submersion. Indeed, restricting to the subset of
metrics which are fixed outside some open subset $U\subset \Sigma$, we
still get a submersion.
\end{prop} 

The proof relies on the following elementary fact:

\begin{lemma}
\label{lemma:LinearAlgebra}
Given a pair of non-zero vectors $v,w\in\R^2$, and a complex structure
$j_0$ on $\R^2$, i.e. an endomorphism of $\R^2$ with $j_0^2 = -\Id$,
there is a one-parameter family $j_t$ of complex structures with
$\DDt\Big|_{t=0} j_t v = w.$
\end{lemma}

\begin{proof}
We can introduce coordinates on $\R^2$ with respect to which $j_0$
takes the form $j_0=\left(\begin{array}{cc}0 & -1 \\ 1 & 0
\end{array}\right)$. For any real numbers $a$, $b$,
let 
$X= \OneHalf \left(\begin{array}{cc} -b & -a \\
-a & b
\end{array}
\right)$. Clearly, $j_t(a,b)=e^{tX}j_0 e^{-tX}$ gives a one-parameter family
of complex structures whose differential at $t=0$ acts by
$$\left(\begin{array}{cc} -a & b \\ b & a
\end{array}\right).$$ The lemma follows.
\end{proof}

\vskip.2cm
\noindent{\bf{Proof of Proposition~\ref{prop:GenericMetrics}}}.
Suppose $h_0$ is some fixed metric, for which $\Phi$ is an
$A_0$-holomorphic section which represents $D$. Then $\Theta_h(D)$ is
defined as follows. Let $a^{0,1}\in\Omega^{0,1}_h$ be a form for which
$$\DBar_{A_0+a^{0,1}}\Phi=0.$$ Note that $a^{0,1}$ depends on the metric
$h$. 
Then, $\Theta_h(D)=\Theta_{\met_0}(D)+ 
[\Proj_{\Harm_h}\Image a^{0,1}]$, where $\Proj_{\Harm_h}$ denotes the
$L^2$ projection map to the space of harmonic one-forms $\Harm_h$ for
the metric $h$. Now if the
metric $h$ differs from $h_0$ only in a region $U\subset \Sigma$ where
$\Phi\neq 0$, then $a^{0,1}$ can be written as follows:
$$a^{0,1}=-\frac{(1+ i J_h) \nabla_{A_0}\Phi}{\Phi}.$$ Here, $J_h$
is (pull-back by) the almost-complex structure induced from the
metric; i.e. it is the Hodge star operator on one-forms.

So, if $h_t$ is a one-parameter family of metrics through $t=0$ (where
$\Phi$ is $A_0$-holomorphic), then the derivative of the theta map is
given by
\begin{eqnarray*}
\DDt_{t=0}\left[\Proj_{\Harm_{h_t}}\Image \left(-\frac{(1+ i J_h)
\nabla_{A_0}\Phi}{\Phi}\right)\right] &=& \left[\Proj_{\Harm_{h_0}} \Image \left( i
\frac{dh_t}{dt}(0) \frac{\nabla_{A_0}\Phi}{\Phi} \right)\right] \\
&=& \left[\Proj_{\Harm_{h_0}} J'(0)\Real \left(
\frac{\nabla_{A_0}\Phi}{\Phi} \right)\right].
\end{eqnarray*}
Now, $b=\Real\left(\frac{\nabla_{A_0}\Phi}{\Phi}\right)$ is a
differential one-form which cannot vanish identically: if it did, that
would mean that $\Phi$ is $A_0$-parallel, since $\Phi$ is
$A_0$-holomorphic. But there are no non-zero, parallel sections of a
bundle of degree $g-1$, unless $g=1$ (which is ruled out by the
hypothesis). 

Even if we choose $J'(0)$ to be supported in the region where $b\neq
0$, we would get a submersion. If this were not the case, we would be
able to find an $h_0$-harmonic one-form $\omega$ which was orthogonal
to the image of that derivative; i.e. for all variations $J'$ of the
metric, we have that $$\langle \Proj_{\Harm_{h_0}}J'\circ b,
\omega\rangle = \langle J'\circ b, \omega\rangle = 0.$$ But if $b,
\omega$ are forms, then one can always find a one-parameter family of
almost-complex structures $J_t$ whose derivative at zero sends $b$ to
a non-negative multiple of $\omega$, in view of
Lemma~\ref{lemma:LinearAlgebra}
\qed


\section{The Dirac operator on the Cylinder}
\label{sec:ALimit}

Let $\{\met_t\}_{t\in[0,1]}$ be a path of metrics on a compact,
oriented two-manifold $\Sigma$ of genus $g$ which are stationary for
all $t<\epsilon$ and $t>1-\epsilon$.  Let $B_t\in\Jac_{\met_t}$ be a
family of $\SpinC$ connections which are stationary for $t<\epsilon$
and $t>1-\epsilon$ as well. Then, these data naturally induce a metric
on $[0,1]\times \Sigma$, together with a connection $A$ on the spinor
bundle.  We introduce a one-parameter family $g_\scale$ of metrics on
$[0,1]\times
\Sigma$, which is to be thought of as ``stretching'' the 
$\R$-factor. Specifically, the metric tensor is given by
$$\metC{\scale}=(\scale dt)^2 + h_t$$
or, equivalently, $\metC{\scale}$ is the pull-back of the metric $dt^2 +
h_{t/\scale}$ over $[0,\scale]\times \Sigma$, by the map
$(t,\sigma)\mapsto (\scale t, \sigma)$.

Since the paths $\{\met_t\}$ and $\{B_t\}$ are $t$-independent near the ends,
they admit natural extensions to $\R\times\Sigma$.

Our aim is to prove the following:

\begin{prop}
\label{prop:StretchOutKernel}
Let $(\{\met_t\}_{t\in[0,1]},B_t)$ be a path of metrics and
connections, with $B_t\in\Jac_{\met_t}$. Suppose moreover that for all
$t\in[0,1]$, the connection $B_t$ misses the $\met_t$ theta
divisor. Then, there is some scale $\scale_0$, so that for all
$\scale\geq \scale_0$, the Dirac operator on $\Cyl$ given the metric
$\metC{\scale}$ has trivial $L^2$ kernel.
\end{prop}

The proposition rests on a ``near-Weitzenb\"ock'' decomposition for
the square of the Dirac operator.

Let $W$ be a spinor bundle over the cylinder (for any of the metrics
$\met_t$). Clearly, the restriction of $W$ to any $t$-slice is a
Clifford bundle over $\Sigma$ with the metric $\met_t$. This allows us
to think of the family of Dirac operators on $\Sigma$ indexed by $t$
(associated to the metric $\met_t$, and connection $B_t$)
as a single operator 
$$\DiracThree\colon W \longrightarrow W,$$ over $\Cyl$. 
Recall that the two-dimensional Dirac operator can be written as 
\begin{equation}
\label{eq:SplitDirac}
\DiracTwo =
\sqrt{2}\left(
\begin{array}{cc}
0 & \DBar_{B_t}^* \\
\DBar_{B_t} & 0
\end{array}
\right)
\end{equation}
with respect to the natural splitting of $W=\SpinBunTwoP\oplus \SpinBunTwoM$.

\begin{lemma}
\label{lemma:NearWeitzenbock}
The Dirac operator for the three-manifold $\R\times \Sigma$ has the
following ``near''-Weitzenb\"ock decomposition
$$\DiracThree^*\DiracThree = -(\frac{1}{\scale}\DDt)^2 +
\DiracTwo^*\DiracTwo + L_\scale \circ \nabla + M_\scale,$$ where the
maps $L_\scale$ and $M_\scale$ are bundle endomorphisms $L_\scale
\colon \Wedge^1 \otimes W \longrightarrow W$ and $M_\scale \colon W
\longrightarrow W$ whose pointwise operator norms go to zero as
$\scale\goesto \infty$; indeed, there are constants $\ConstOne$ and $\ConstTwo$,
so that for all $\scale$, we have
\begin{eqnarray*}
\|L_{\scale}\|\leq
\frac{1}{\scale} \ConstOne,
&{\text{and}}&\|M_{\scale}\|=\frac{1}{\scale}\ConstTwo.
\end{eqnarray*}
(The bounds here are pointwise bounds on sections of endomorphism bundles.)
\end{lemma}

In the proof of the lemma, we find it convenient to use a reducible
connection on $T(\Cyl)$, which is independent of $\scale$, defined as
follows. Consider the natural orthogonal splitting (which is valid for
all $\scale$) $$\Wedge^1 \Cyl = \R dt \oplus \Proj_{\Sigma}^*
\Wedge^1.$$ The family of Levi-Civita connections on $\Sigma$ which
arises from the family of metrics $\met_t$ can be viewed as a single
connection on the $dt^\perp$ summand. Now, let $\TConn$ denote the
(reducible) connection which is the connect sum of that connection
with the connection on the trivial line bundle $\R dt$ for which
$dt$ is covariantly constant. If $W_0$ is a spinor bundle over $\Sigma$,
consider the bundle $W=\Proj_\Sigma^* W_0$ over $Y$. $W$ already has a
Clifford action of $dt^\perp$. For each $\scale$, we can complete this
action $$\rho_\scale \colon TY\otimes W \longrightarrow W,$$ by
defining $$\rho_\scale (\scale dt) = \left(\begin{array}{cc} i & 0 \\
0 & -i
\end{array}\right).$$ 

\begin{proof} 
We calculate the (connection form) difference between the Levi-Civita
connection and $\TConn$ over $T(\Cyl)$. Suppose for the moment that
$\scale=1$.
Let $\theta^1, \theta^2$ be a
(time dependent) moving coframe for $\Sigma$.
There are functions $w,x,y,z$ over $\Cyl$ for which 
\begin{eqnarray*}
d\theta^1&=& wdt\wedge \theta^1 + x dt\wedge\theta^2 + \omega^1_2
\wedge \theta^2  \\
d\theta^2&=& ydt\wedge \theta^1 + z dt\wedge\theta^2 -
\omega^1_2\wedge \theta^1,
\end{eqnarray*}
where $\omega^1_2$ is the connection matrix for the metric $\met_t$
over $\Sigma$.
Thus, by the Cartan formalism, the connection matrix with respect to the
coframe $(dt,\theta^1,\theta^2)$ is given by
\begin{equation}
\label{eq:ConnMatrix}
\left(
\begin{array}{ccc}
0 & -w\theta^1 - \left(\frac{x+y}{2}\right)\theta^2 & -z \theta^2 -
\left(\frac{x+y}{2}\right)\theta^1 \\
w \theta^1 + \left(\frac{x+y}{2}\right)\theta^2 & 0 &
\left(\frac{y-x}{2}\right)dt + \omega^1_2\\
z\theta^2 +\left(\frac{x+y}{2}\right)\theta^1 &
\left(\frac{x-y}{2}\right)dt +\omega^2_1 & 0 
\end{array}
\right).
\end{equation}

For general $\scale$, $(\scale dt, \theta^1, \theta^2)$ form a moving
coframe, and the above calculation shows that the difference form
$$\TConn_\scale-\LConn_\scale = \frac{1}{\scale} \Xi_\scale,$$
where $\Xi_\scale$ is a one-form whose $g_\scale$ length is
independent of $\scale$. In fact, by glancing at the connection
matrix, we have that 
$\PDiracThree_\scale-\DiracThree_\scale = \xi,$
for some endomorphism $\xi$ which is independent of $\scale$. 
It follows that the square can be written
$$\DiracThree\circ \DiracThree = - \frac{1}{\scale^2}\DDt\circ \DDt +
\DiracTwo^*\circ \DiracTwo + \{\rho_\scale(\scale dt)\nabla_{{\scale
dt}^\flat}, \DiracTwo\} + 
\frac{1}{\scale}(\xi^* \PDiracThree  + \PDiracThree \xi).$$

The anticommutator term can be expressed in the local coframe:
\begin{eqnarray*}
\{\rho_\scale(\scale dt)\nabla_{{\scale
dt}^\flat}, \DiracTwo\} &=&
\{\rho_\scale(\scale dt)\nabla_{{\scale
dt}^\flat}, \theta^1\nabla_{\theta^1}+ \theta^2\nabla_{\theta^2}\} \\
&=&
\rho_{\scale}(dt\wedge\theta^1)[\nabla_{\scale dt},\nabla_{\theta^1}] 
+\rho_{\scale}(dt\wedge\theta^2)[\nabla_{\scale
dt},\nabla_{\theta^2}] \\
&=&
\rho_{\scale}(\Proj_{(\Wedge^2_\Sigma)}^\perp\circ F_\scale)
+ \frac{1}{\scale}\left(\rho_{\scale}(\scale dt
\wedge\theta^1)\nabla_{[\DDt,e_1]} + \rho_{\scale}(\scale dt
\wedge\theta^2)\nabla_{[\DDt,e_2]}\right).
\end{eqnarray*}
The two-form $F_\scale$ is independent of $\scale$ (only its Clifford
action depends on $\scale$); and the Clifford action of a fixed form
in
$dt\wedge \Wedge^1_\Sigma\subset TY$ scales like $1/\scale$.

\end{proof}

In view of Equation~\eqref{eq:SplitDirac}, the hypothesis that $B_t$
always misses the theta divisor is equivalent to the statement that
there is a non-zero lower bound $\delta_0>0$ on the square of
the eigenvalues of the two-dimensional Dirac operator
$\DiracTwo$; in particular,
$$\int_{\{t\}\times \Sigma}\langle \DiracTwo \Psi,
\DiracTwo\Psi
\rangle \geq \delta_0 \int_{\{t\}\times\Sigma}\langle
\Psi,\Psi\rangle.$$ This, together with the formula from
Lemma~\ref{lemma:NearWeitzenbock}, gives us a differential inequality
for harmonic spinors on the cylinder.

\vskip.3cm
\noindent{\bf{Proof of Proposition~\ref{prop:StretchOutKernel}.}}
Let $\Phi$ be a kernel element of $\DiracThree$. Then,
\begin{eqnarray*}
0 &=&  \DiracThree\circ\DiracThree~\Phi \\
&=& -\frac{1}{\scale^2} \frac{\partial^2}{\partial t^2} \Phi + \DiracTwo^*\DiracTwo \Phi + L_\scale \circ \nabla \Phi +
M_\scale \Phi.
\end{eqnarray*}
Taking pointwise inner-product with $\Phi$ (using the metric on $W$),
and integrating out over $\{t\}\times \Sigma$, we see that
\begin{eqnarray}
0 &=& \langle -\frac{1}{\scale^2}\DDtSq \Phi , \Phi\rangle + \langle
\DiracTwo\Phi,\DiracTwo\Phi\rangle + \langle L_\scale
\circ \nabla\Phi,\Phi \rangle +
\langle M_\scale \Phi,\Phi\rangle \nonumber\\
&\geq& - \frac{1}{2\scale^2}\DDt^2 \langle \Phi,\Phi\rangle + \left(\delta_0 -
\frac{1}{\scale}\ConstTwo\right) \langle
\Phi,\Phi \rangle - \langle L_\scale\circ \nabla\Phi,\Phi \rangle.
\label{eq:DiffIneqTakeOne}
\end{eqnarray}
Now, we have:
\begin{eqnarray}
\left|\int_{\{t\}\times \Sigma}
\langle L_\scale \circ\nabla\Phi,\Phi \rangle \right|
&\leq&
\frac{1}{\scale} \ConstOne \left(\int_{\{t\}\times \Sigma} |\nabla\Phi|^2\right)^\OneHalf
\left(\int_{\{t\}\times \Sigma}|\Phi|^2\right)^\OneHalf
\label{eq:BoundLTermTakeOne}
\end{eqnarray}
But, we have, by integration-by-parts and the usual Weitzenb\"ock formula: 
\begin{eqnarray*}
\langle \nabla\Phi,\nabla\Phi\rangle &=& \frac{1}{2\scale^2} \DDtSq \langle \Phi,
\Phi\rangle + \langle \nabla^*\nabla
\Phi,\Phi\rangle  \\
&=& \frac{1}{2\scale^2} \DDtSq \langle \Phi,
\Phi\rangle + \langle (s_\scale + \rho_\scale(F_A))
\Phi,\Phi\rangle  \\
&\leq & \frac{1}{2\scale^2} \DDtSq \langle \Phi,
\Phi\rangle + C \langle \Phi,\Phi\rangle,
\end{eqnarray*}
for some non-negative constant $C$ independent of $\scale$
(here, $s_{\mu}$ is the scalar curvature of the stretched manifold; 
this curvature stays bounded as $\mu\goesto\infty$);
so, substituting this into Inequality~\eqref{eq:BoundLTermTakeOne},
then plugging back into Equation~\eqref{eq:DiffIneqTakeOne}, we get
the following differential inequality for the function 
$\|\Phi\|^{2}$ (thought of as a function of $t$, given by
$t \mapsto \int_{Y}|\Phi(t,y)|^{2}dy$):
\begin{eqnarray*}
0 &\geq& - \frac{1}{2\scale^2}\DDt^2 \|\Phi\|^2 + \left(\delta_0 -
\frac{1}{\scale}\ConstTwo\right) \|\Phi\|^2 -
\frac{1}{\scale}\ConstOne\left(\frac{1}{2\scale^2}\DDtSq{\|\Phi\|^2} +
C\|\Phi\|^2\right)^\OneHalf \|\Phi\|.
\end{eqnarray*}
By rearranging the above inequality we get for sufficiently large $\scale$
an inequality of the form $$0\geq -\frac{1}{2\scale^2}\DDt^2
\|\Phi\|^2 + {\delta}\|\Phi\|^2,$$ where $\delta>0$ is some constant
(which we can make arbitrarily close to $\delta_0$ by choosing
$\scale$ to be sufficiently large). This in turn shows that
$\|\Phi\|^2$ has no interior maxima.  In particular, if $\Phi$ is
integrable, it follows that $\Phi\equiv 0$.
\qed
\vskip.3cm

We now spell out:

\begin{prop}
\label{prop:SFCylBig}
Let $\met_t$, $\met_t'$ be two generic allowable paths of metrics over
$\Sigma$, and fix an allowable homotopy $H(u,t)$ between them.
Suppose that $B(u,t)$ is a two-parameter family of connections in 
$\Jac_{H(u,t)}$
which is transverse to $\Theta_H$.
Then, there is some $\scale_0\geq 0$ such that  
for all sufficiently large $\scale,\scale'\geq \scale_0$, 
the spectral flow
between the two induced Dirac operators on $\Cyl$ is given by
\begin{eqnarray*}
\lefteqn{\SF\left((\mu dt)^2+\met_t, (\mu'dt)^2 + \met_t'\right)} 
\\
&=& \#\left\{(D,u,t)\in\Sym^{g-1}(\Sigma)\times[0,1]\times[0,1]\big|
\Theta_{H(u,t)}(D)=B(u,t)\right\}.
\label{eq:IntNumSF}
\end{eqnarray*}
\end{prop}

\begin{proof}
By Proposition~\ref{prop:StretchOutKernel}, the spectral flow
localizes to the region where $B(u,t)$ meets the $H(u,t)$-theta divisor.
By homotopy
invariance of both quantities, then, it follows that the spectral flow
must be some multiple of the above intersection number. The factor is
then calculated in a model case, as in
Proposition~\ref{Theta:prop:SFInterp} of~\cite{Theta} (note that 
in that proposition, the ``spectral flow'' refers to real spectral 
flow, hence the difference in factors of $2$).
\end{proof}


\section{$\Ntheta$ and the Casson invariant}
\label{sec:Casson}

In this section, we prove Theorem~\ref{thm:Casson}.  At the heart of
this computation is a surgery formula.

We focus presently on the case where $Y$ is an integral homology
three-sphere. In this case, there is only one $\SpinC$ structure, and
we denote its invariant by $\Ntheta(Y)$. Note that $\Ntheta(Y)$ is
an integer, since in this case, the APS correction term $\Corr(Y)$ is
an integer for all metrics.

Let $Y$ as above and fix a knot $K\subset Y$. Let $Y_{p/q}$ be the
manifold obtained from $Y$ by $(p/q)$-surgery along $K$, so that $Y_0$
is an integral homology $S^1\times S^2$, and $Y_1$ is another integral
homology three-sphere. Note that, since $b_1(Y_0)=1$, the theory
from~\cite{Theta} applies, to give us an integer-valued invariant
$\theta(Y_0,\spinc)$ for each $\SpinC$ structure
$\spinc\in\SpinC(Y_0)$.

Our first result is the following surgery formula for $\Ntheta$.

\begin{theorem}
\label{thm:SurgeryFormula}
$$\Ntheta(Y_1)-\Ntheta(Y)=\sum_{\spinc\in\SpinC(Y_0)}\theta(Y_0,\spinc).$$
\end{theorem}

This theorem, together with Theorem~\ref{Theta:thm:CalculateZ}
from~\cite{Theta}, which relates $\theta(Y_0)$ with the Alexander
polynomial of $Y_0$, gives the following:

\begin{theorem}
Let $A=a_0+\sum_{i=1}^k a_i(T^i+T^{-i})$ be the symmetrized Alexander
polynomial of $Y_0$, normalized so that
$$A(1)=1.$$ Then,
$$\Ntheta(Y_1)-\Ntheta(Y)=\sum_{j=1}^{\infty} j^2 a_j.$$
\end{theorem}

\begin{proof}
We recall from~\cite{Theta} that $\theta(Y_0,\spinc)=\sum_{j=1}^\infty
j a_{|i|+j}$, where the $\SpinC$ structure $\spinc$ corresponds to the
integer $i$ under the identification $\SpinC(Y_0)\cong H^2(Y;\Z)\cong
\Z$ which sends the spin structure to $0$. 
\end{proof}

\begin{remark}
It is an easy consequence of this theorem that:
$$\Ntheta(Y_{1/n})-\Ntheta(Y)=n\sum_{j=1}^{\infty} j^2 a_j.$$
\end{remark}

\begin{cor}
Let $Y$ be an oriented homology three-sphere. Then, $2\Ntheta(Y)$ is
equal to the Casson invariant $\lambda(Y)$.
\end{cor}

\begin{proof}
Every integral homology three-sphere can be obtained from $S^3$ by a
sequence of $\pm 1$ surgeries. So, it follows that $\NInv$ is uniquely
determined by its surgery formula and its value on $S^3$.  It is easy
to see that $\NInv(S^3)=0$: fix a genus one Heegaard decomposition of
$S^3$, and let $\met_t$ be a constant family of metrics on the torus,
which is clearly an allowable path, and indeed $\Inv(S^3)=0$ for this
path.  Furthermore, the correction terms cancel each other by
Proposition~\ref{prop:CompareCorrections}.  

Since Casson's invariant $\lambda$ satisfies the same surgery formula
(see~\cite{Casson}), and $\lambda(S^3)=0$ as well, we get the result.
\end{proof}

To prove Theorem~\ref{thm:SurgeryFormula}, we find a Heegaard
decomposition $Y=U_0\cup_{\Sigma} U_1$ for which the knot $K$ is dual to
the last attaching disk for $U_1$, i.e. it intersects the last
attaching disk transversally in one point, and is disjoint from the
other ones.  To see that this can be arranged, start with a Morse
function on $Y-\nbd{K}$ with $1$ zero-handle, $g$ one-handles, and
$g-1$ two-handles. This gives us a handlebody $U_0$ and $g-1$
two-handles, whose attaching circles we denote by
$\{\beta_1,...,\beta_{g-1}\}$. Completing the Morse function over
$\nbd{K}$, we have described the second handlebody $U_1$ so that $K$
is dual to the final attaching disk. Heegaard decompositions for the
various surgeries $Y_{p/q}$ are given by surgeries in $U_1$: thus,
these can be described by fixing one surface $\Sigma$, one complete
set of attaching circles $\{\alpha_1,...,\alpha_g\}$, a $g-1$-tuple of
attaching circles $\{\beta_1,...,\beta_{g-1}\}$, and allowing the
final attaching circle $\beta_g$ to vary.

We would like to compare the $\theta$-invariants for the various
surgeries. To that end, we fix a path $\met_t$ of metrics, so that
$\met_0$ is $U_0$-allowable and $\met_1$ is $U_1$ allowable for any
choice of $\beta_g$. Such a metric $\met_1$ can be found, thanks to
Lemma~\ref{Theta:lemma:MissThetaHS} of~\cite{Theta}, which shows that
any metric which is sufficiently stretched out normal to the
$\{\beta_1,...,\beta_{g-1}\}$ is $U_1$-allowable.

We define a ``moduli space'' belonging to the knot complement and the
path of metrics.
$$\ModSp_{\met_t}(Y-K)=\{(D,s,t)\in\Sym^{g-1}\times[0,1]\times[0,1]\big|
s\leq t, 
\Theta_{\met_s}(D)\in L(U_0),
\Theta_{\met_t}(D)\in\Lambda(U_1)\},$$
where $\Lambda(U_1)\subset \Jac$ is the set of connections $B\in\Jac$
so that $\Hol_{\beta_i}(B)=0$ for $i=1,...,g-1$.  There is a 
map $$\rho\colon \ModSp_{\met_t}(Y-K)\longrightarrow
\Torus{2},$$
which is analogous to a boundary value map, given by measuring
holonomy around the meridian $m$ and the longitude $\ell$ of $K$,
i.e.  $$\rho(D,s,t)=\Hol_{m\times\ell}(\Theta_{\met_t}(D)),$$
normalized so that the point $0\times 0 \in \Torus{2}$ corresponds to
the spin structure which extends over $Y$.  With these conventions,
then, any ``reducible'' (in the sense of 
Definition~\ref{def:Reducible})
on $Y-K$ can be restricted to
the boundary; its holonomies will lie in the circle $S^1\times\{0\}$.

Note that $\Lambda(U_1)\cap \Theta_{\met_1}(\Sym^{g-1}(\Sigma))$ is
not necessarily empty. However, for any $\epsilon>0$, there is a $T$
so that if $\met_1$ is stretched out at least $T$ normal to the
$\{\beta_1,...,\beta_{g-1}\}$, then the holonomy of any point in
$\Lambda(U_1)\cap \Theta_{\met_1}(\Sym^{g-1}(\Sigma))$ around
$m\times\ell$ lies in an $\epsilon$ neighborhood of
$\OneHalf\times\OneHalf$. This follows from
Lemma~\ref{Theta:lemma:MissThetaHS} of~\cite{Theta}.

This moduli space has the following important properties:

\begin{prop}
\label{prop:StructureModSpace}
For any $\epsilon>0$, if $\met_1$ is sufficiently stretched out normal
to the attaching circles $\{\beta_1,...,\beta_{g-1}\}$, then the
moduli space of the knot complement $\ModSp_{\met_t}(Y-K)$ is
generically a compact, smooth, one-dimensional manifold with two types
of boundary components corresponding to $t=1$ and $s=t$.  Furthermore,
the $t=1$ boundary maps under $\rho$ into an $\epsilon$-neighborhood
of $\OneHalf\times\OneHalf$, and those with $s=t$ map under $\rho$ to
$S^1\times 0$.
\end{prop}

\begin{proof}
Smoothness follows from the generic metrics statement
(Proposition~\ref{prop:GenericMetrics}). There is no $s=0$ boundary
since the metric $\met_0$ is $U_0$-allowable. The $t=1$ boundary lies in
$\Lambda(U_1)\cap \Theta_{\met_1}(\Sym^{g-1}(\Sigma))$, which maps near
$\OneHalf\times\OneHalf$, as above. The $s=t$ boundary corresponds to
the intersection of the theta divisor with $\Lag(U_0)\cap
\Lambda(U_1)$, which in turn maps to $S^1\times 0$. Indeed, this
latter circle corresponds to the circle of reducibles
$Y-K$: 
the point $0\times 0$ corresponds to the spin structure
which extends over $Y$, while the point $\OneHalf\times 0$ corresponds
to the spin structure on $Y_1$.
\end{proof}

As we shall see presently, the $\theta$-invariants for the surgered
manifolds $Y$, $Y_1$, and $Y_0$ are related by a spectral flow term,
defined as follows. Let $(\scale dt)^2+\met_t$ denote the metric on
the cylinder $\R\times\Sigma$ induced by the family $\met_t$. By
restriction, the spin structures $\spinc_0$ and $\spinc_1$ on $Y$ and
$Y_1$ respectively induce spin structures on $Y-K$.  These can be
viewed as different $\SpinC$-connections on the same $\SpinC$
structure on $Y-K$, which we will connect by a path
$\{A_t\}$ of reducible connections. For definiteness,
we choose a path whose holonomies around the meridian are monotone
increasing from $0$ to $\OneHalf$.  By restricting to the cylinder
$\R\times\Sigma$, we get a path of reducibles $\{A_t\}$, 
and hence a corresponding path
of Dirac operators $\Dirac_{A_t|\R\times\Sigma}$.  As we have seen
(Proposition~\ref{prop:SFCyl}), if the metric on $\R\times\Sigma$ is
sufficiently stretched out, then this spectral flow between the Dirac
operators is independent of the scale $\scale$. We denote this
spectral flow by $\SF_{\R\times\Sigma}(\spinc_0,\spinc_1)$.  With this
spectral flow defined, we turn our attention to the following key step
towards establishing Theorem~\ref{thm:SurgeryFormula}:

\begin{prop}
\label{prop:SurgerTheta}
Fix a generic path of metrics $\met_t$ so that $\rho(\ModSp_{\met_t}(Y-K))$
misses the spin structures $0\times 0$ and $\OneHalf\times 0$.  The
theta invariant satisfies:
\begin{equation}
\label{eq:SurgerTheta}
\theta_{\met_t}(Y_1)-\theta_{\met_t}(Y)=\left(\sum_{\spinc\in\SpinC(Y_0)}\theta(Y_0,\spinc)\right)
+ \SF_{\R\times\Sigma}(\spinc_0,\spinc_1).
\end{equation}
\end{prop}

\begin{proof}
Fix a real number $\delta$ with $0<\delta<\OneHalf-\epsilon$.  Fix
curves in $\Torus{2}$: $\gamma= 0\times S^1$, $\gamma_0=S^1\times
\{\delta\}$, and $\gamma_1=\{(s,s+\OneHalf)\}$. Note that the restriction of $\spinc_0$ and $\spinc_1$ 
to the torus gives the spin structures corresponding to $0\times 0$
and $\OneHalf\times 0$ respectively.  Since the moduli space misses
these spin structures, it follows that $\theta_{\met_t}(Y)$ and
$\theta_{\met_t}(Y_1)$ are well-defined, and indeed by the definition
of the $\theta$-invariant, we have that
$\#\rho^{-1}(\gamma)=\theta_{\met_t}(Y)$ and
$\#\rho^{-1}(\gamma_1)=\theta_{\met_t}(Y_1)$; similarly,
$$\#\rho^{-1}(\gamma_0)=
\left(\sum_{\spinc\in\SpinC(Y_0)}\theta(Y_0,\spinc)\right).$$

Consider the oriented subset $\Dom\subset \Torus{2}$ which does not
contain $\OneHalf\times\OneHalf$, and whose boundary is $\gamma_0 +
\gamma - \gamma_1$. Then, by transversality, 
\begin{equation}
\label{eq:BoundModSp}
\partial \left(\rho^{-1}(\Dom)\right)=
\rho^{-1}(\gamma_0)+\rho^{-1}(\gamma)-\rho^{-1}(\gamma_1)+
\left(\partial\ModSp_{\met_t}(Y-K)\right)\cap \rho^{-1}(\Dom).
\end{equation}
According to Proposition~\ref{prop:StructureModSpace}, and our choice
of $\Dom$, 
$\left(\partial\ModSp_{\met_t}(Y-K)\right)\cap \rho^{-1}(\Dom)$ consists of
boundary components where $s=t$. In fact,
Proposition~\ref{prop:SFCylBig} shows that 
$$\#\left(\partial\ModSp_{\met_t}(Y-K)\right)\cap
\rho^{-1}(\Dom)=\SF_{\R\times\Sigma}(\spinc_0,\spinc_1),$$
so that counting points in Equation~\eqref{eq:BoundModSp}, we obtain
Equation~\eqref{eq:SurgerTheta}.
\end{proof}

We can understand the spectral flow on the cylinder in terms of data
on the knot complement, thanks to the splitting principle for spectral
flow. Specifically, we can connect $\spinc_0$ and $\spinc_1$ through
reducibles on the knot complement, endowed with a metric which is
$\Umet_0$ on $U_0$, $\met_t$ on the cylinder, and a fixed metric on
$U_1-\nbd{K}$ which is product-like near both the $\Sigma$-boundary
(where it is isometric to a collar around $\met_1$), and the torus
boundary of $\nbd{K}$. (In fact, for concreteness, we assume that it
has the form $S^1\times S^1$ near the knot, where we take the product
of longitude and meridian.)  According to the splitting principle,
then,
\begin{equation}
\label{eq:SplitKnotK}
\SF_{Y-K}(\spinc_0,\spinc_1)=\SF_{U_0}(\spinc_0,\spinc_1)+\SF_{\R\times\Sigma}(\spinc_0,\spinc_1)
+ \SF_{U_1-K}(\spinc_0,\spinc_1).
\end{equation}

Strictly speaking, in order to achieve the apropriate transversality,
so that the Dirac operator at the endpoints $\spinc_0$ and $\spinc_1$
have no kernel, we have relaxed the reducibility hypothesis, to include
perturbations of reducibles by one-forms which are compactly supported
in the two handlebodies (see the proof of Lemma~\ref{lemma:UDirac}).
Thus, for example, Equation~\eqref{eq:SplitKnotK} should be
interpreted as follows: fix a pair $a_0$ and $a_1$ of one-forms which
are compactly supported inside $U_0$ and $U_1-K$ respectively.
$\SF_{Y-K}(\spinc_0,\spinc_1)$ is the spectral flow of the family $A_t
+ a_0 + a_1$, where $A_t$ is the family of reducibles connecting
$\spinc_0$ to $\spinc_1$, and the spectral flows over $U_0$ and
$U_1-K$ are actually the restrictions of $A_t+a_0+a_1$ to these two
submanifolds.  In the interest of clarity, we suppress these
perturbations from the subsequent discussion: when we discuss to
``reducibles'', we mean, in fact, connections obtained by
perturbing reducibles in this manner. (Note also that these
perturbations leave the holonomies around attaching circles
unchanged.)

To relate the spectral flow on the knot complement with data on closed
manifolds, we have the following result, which could be called a
surgery formula for $\Corr$ (compare~\cite{RubermanMeyerhoff}, see
also~\cite{Lim}):

\begin{lemma}
\label{lemma:SurgerCorr}
Fix a cylindrical-end metric $\Kmet$ on $Y-K$, and let $\metU_0$,
$\metU_1$ be metrics on the genus one handlebody $U$ which agree with
$\Kmet$ along the boundary.  Then, for all sufficiently large $T$,
$$\Corr_{Y_1}(\Kmet\cup_{T}\metU_1)-\Corr_{Y}(\Kmet\cup_{T}\metU_0)=
\SF_{Y-K}(\spinc_0,\spinc_1) + \CorrU(k_1, \spinc_1)-\CorrU(k_0, 
\spinc_0).$$
\end{lemma}

\begin{proof}
Recall that there is a natural cobordism $W$ from $Y$ to $Y_1$.  This
cobordism is obtained from $Y\times [0,1]$, and then attaching a
two-handle with $+1$ framing.  We have a natural inclusion
$$(Y-K)\times [0,1]\subset W,$$ which maps $(Y-K)\times \{0,1\}$ to
the boundary of $W$ (i.e. taking $(Y-K)\times\{0\}$ to the knot
complement, as a subset of $Y$, and $(Y-K)\times\{1\}$ to the
corresponding subset of $Y_1$). Fix a metric on $W$ whose restriction
to the region
$(Y-K)\times[0,1]$ is a product metric, its restriction to the $Y$
boundary agrees with $\Kmet\cup_T k_0$, and its restriction to $Y_1$
agrees with $\Kmet\cup_T k_1$.

Consider the $\SpinC$ structure $\spincX$ on $W$ whose first Chern
class generates $H^2(W;\Z)$. Endow $W$ with a $\SpinC$ connection in
$\spincX$ whose restriction to $(Y-K)\times [0,1]$ is a path of
reducibles, and whose first Chern form is compactly supported away
from the boundary (and hence it interpolates between $\spinc_0$ and
$\spinc_1$). By definition,
\begin{equation}
\label{eq:DifferenceCorrY}
\Corr_{Y_{1}}-\Corr_{Y} = \ind \Dirac_{A} -
\left(\frac{c_1(\spincX)^2-\sigma(W)}{8}\right)=\ind\Dirac_{A}.
\end{equation}

The index is calculated by an excision principle. Let
$Z=W-(Y-K)\times[0,1]$.  According to~\cite{MullerCorners}, the Dirac
operator on both $Z$ and $(Y-K)\times[0,1]$ is a Fredholm operator,
since the Dirac operator has no kernels on the ``corners'' $S^1\times
S^1\times\{0,1\}$ and the various boundaries $Y-K$ (coupled to
$\spinc_0$ and $\spinc_1$ respectively) and $[0,1]\times \Sigma_1$
(coupled to a sufficiently slowly-moving one-parameter family of
connections on $\Sigma_1$ which bound -- see
Proposition~\ref{prop:StretchOutKernel}).  Thus, the index splits as
$$\ind\Dirac_A = \ind {\Dirac_A}|_{(Y-K)\times[0,1]} +
\ind{\Dirac_A}|_{Z}.$$
The first term is the spectral flow $\SF_{Y-K}(\spinc_0,\spinc_1)$. To
understand the second term, we replace the knot complement $Y-K$ by a
much simpler knot-complement $\CDisk\times S^1$, endowed with a family
of connections whose holonomy around the $S^1$ factor goes from $0$ to
$\OneHalf$. Note that the manifold $W_0=\left(\CDisk\times
S^1\right)\cup_{S^1\times S^1} Z$ is a cobordism from $S^3$ to $S^3$:
indeed, it is $\CP{2}$ punctured at two points. Moreover, the connection $A|_Z$ extends over $W_0$, to give a connection $A_0$. Now, by the same splitting
principle, 
$$\ind\Dirac_{A_0} = \ind {\Dirac_{A_0}}|_{(\CDisk\times S^1)\times[0,1]} +
\ind{\Dirac_{A_0}}|_{Z}.$$
The first term on the right hand side is a spectral flow for the Dirac
operator through flat connections on the manifold $\CDisk\times S^1$,
which has non-negative sectional curvatures; thus, the index vanishes. The
second term on the right is, of course, the same as
$\ind{\Dirac_A}|_{Z}$; thus, we have that
\begin{equation}
\label{eq:Index}
\ind\Dirac_{A_0}=\ind {\Dirac_A}|_Z
\end{equation}
Finally, by the same reasoning which gave
Equation~\eqref{eq:DifferenceCorrY}, we have that
\begin{equation}
\label{eq:DiffCorrS}
\ind{\Dirac_{A_0}}=\Corr((\CDisk\times S^1)\cup
(U,k_0))-\Corr((\CDisk\times S^1)\cup (U,k_1)),
\end{equation}
which is a difference of the APS invariant for the three-sphere
$S^3$. Now, according to Proposition~\ref{prop:CompareCorrections}
(together with the fact that $\CorrU(\CDisk\times S^1)\equiv 0$ for
any $\SpinC$ connections with traceless curvature), we have that
\begin{equation}
\label{eq:SplitCorrS}
\Corr((\CDisk\times S^1)\cup U_1)
-
\Corr((\CDisk\times S^1)\cup U_0)
=
\CorrU(U_0)-\CorrU(U_1).
\end{equation}
Combining Equations~\eqref{eq:DifferenceCorrY}-\eqref{eq:SplitCorrS},
we have established the lemma.
\end{proof}

The surgery formula for $\NInv$ now is an easy consequence of the
surgery formula for $\theta$ and $\Corr$:

\vskip.2cm
\noindent{\bf Proof of Theorem~\ref{thm:SurgeryFormula}.}
First we claim that the chambered properties of $\CorrU$ and $\Corr$,
the splitting formula (the version stated in
Equation~\ref{eq:SplitKnotK}), and Lemma~\ref{lemma:SurgerCorr}, we
see that 
$$
\Big(\CorrY(Y_1,\spinc_1)-\CorrU(U_0,\spinc_1)-\CorrU(U_1',\spinc_1)\Big) 
-
\Big(\CorrY(Y,\spinc_0)-\CorrU(U_0,\spinc_0)-\CorrU(U_1,\spinc_0)\Big)
= \SF_{\R\times
\Sigma}(\spinc_0,\spinc_1).
$$ It follows that
$$\NInv(Y_1)-\NInv(Y)=\Inv_{\met_t}(Y_1)-\Inv_{\met_t}(Y) -
\SF_{\R\times \Sigma}(\spinc_0,\spinc_1).$$
Proposition~\ref{prop:SurgerTheta} then finishes the proof.
\qed


\section{$\NInv$ and the Casson-Walker invariant}
\label{sec:Walker}

Let $X$ be an oriented three-manifold with a torus boundary and
$H_1(X;\R)\cong \R$.  The map $H_1(\partial X;\Z)
\longrightarrow H_1(X;\Z)$ has one-dimensional  kernel. Let
$\ell'$ denote a generator for the kernel, $d(X)>0$ denote its
divisibility, and let $\ell$ be the element $\ell'/d$. We call $\ell$
the {\em longitude}.

Fix a homology class $m\in H_1(\partial X)$ with $m\cdot
\ell=1$.  For a pair of relatively prime integers $(p,q)$, the
manifold $Y_{p/q}$ is obtained from $X$ by attaching a
$S^1\times \CDisk$ with $\partial \CDisk= p m + q \ell$, and let
$Y=Y_{1/0}$. Note that in general $Y_{p/q}$ depends on a choice of
$m$, but $Y_{0}=Y_{0/1}$ does not. Note also that $Y_{0}$ is a
rational homology $S^1\times S^2$, while all the other $Y_{p/q}$ are
rational homology spheres.

There is a short exact sequence
\begin{equation}
\label{eq:SESYnot}
\begin{CD}
0@>>> \Z @>>> \SpinC(Y_0) @>>> \SpinC(X) @>>> 0
\end{CD},
\end{equation}
by which we mean that the subgroup $\Z\subset H^2(Y_0;\Z)$ generated
by the Poincar\'e dual to $m$ (viewed as a subset of $Y_0$) acts
freely on $\SpinC(Y_0)$, and its quotient is naturally identified
(under restriction to $X\subset Y_0$) with $\SpinC(X)$.

Thus, each $\SpinC$ structure $\spinc$ on $X$ has a natural level
$y=y(\spinc)\in\Zmod{d}$ defined as follows.  Let $\spincb$ be any
$\SpinC$ structure on $Y_0$ whose restriction is $\spinca$, and
consider its image in $$\tSpinC(Y_0)/\Z(\PD[m])\cong \Zmod{d},$$
where $\tSpinC(Y_0)$ is the group of $\SpinC$ structures modulo the
action of the torsion subgroup of $H^2(Y_0;\Z)$. 

Furthermore, for any of the $Y_{p/q}$, the map $\SpinC(Y_{p/q})$ to
$\SpinC(X)$ is surjective, and its fibers consist of orbits by a
cyclic group generated by the Poincar\'e dual to the knot which is the
core of the complement $Y_{p/q}-X$ (for $Y=Y_{1/0}$, this fiber has
order $d=d(X)$). For a fixed $\SpinC$ structure $\spinc$ on $X$, let
$\SpinC(Y_{p/q};\spinc)$ denote the set of $\SpinC$ structures
$\spincb\in\SpinC(Y_{p/q})$ whose restriction to $X$ is $\spinc$.

Our main result in this section is:

\begin{theorem}
\label{thm:PQSurgeryFormula}
For integers $p,q,d,y$ with $p$ and $q$ relatively prime, $d>0$ and
$0\leq y< d$, there is quantity $\epsilon(p,q,d,y)\in\Q$ with the
following property.  Let $X$ be an oriented rational homology
$S^1\times \CDisk$, with divisibility $d(X)=d$, and choose $m$, $\ell$
as above.  Fixing any $\SpinC$ structure $\spinc$ over $X$ with level
$y(\spinc)=y$, we have the relation:
$$\left(\sum_{\spincb\in\SpinC(Y_{p/q};\spinc)}\NInv_{Y_{p/q}}(\spincb)\right)
= p \left(\sum_{\spincc\in\SpinC(Y;\spinc)}\NInv_{Y}(\spincc)\right) +
q
\left(\sum_{\spincd\in\SpinC(Y_0;\spinc)}\Inv_{Y_0}(\spincd)\right) 
+\epsilon(p,q,d,y).$$
\end{theorem}

\begin{cor}
\label{cor:PQSurgeryFormulaSum}
For $X$ as above,
$$\left(\sum_{\spincb\in\SpinC(Y_{p/q})}\NInv_{Y_{p/q}}(\spincb)\right)
= p \left(\sum_{\spinc\in\SpinC(Y)}\NInv_{Y}(\spinc)\right)+
q(\sum_{i=1}^\infty a_i i^2) + |\Tors H_1(X;\Z)| \epsilon(p,q,d),$$
where $d=d(X)$, $a_i$ are the coefficients of the
symmetrized Alexander polynomial of $Y_0$, normalized so that
$$
A(1)=|\Tors H^2(Y_0;\Z)|,
$$
and 
$$\epsilon(p,q,d)=\sum_{y=0}^{d-1}\frac{\epsilon(p,q,d,y)}{d}.$$
\end{cor}

\begin{proof}
This follows from the surgery formula, and the calculation of 
$\theta(Y_{0})$ from~\cite{Theta}, according to which
$$ \sum_{\spinc\in\SpinC(Y_0)}\Inv_{Y_0}(\spinc) =
\sum_{i=1}^\infty a_i i^2.$$
\end{proof}

\begin{remark}
\label{rmk:Surgeries}
Let $K\subset Z$ be a knot in any rational homology three-sphere. In
the above notation, $Z=Y_{p/q}$ for $X=Z-\nbd{K}$ and some choice of
$m$, $p$, and $q$. Clearly, any non-zero surgery on $K$ would give
another rational three-sphere of the form $Z'=Y_{p'/q'}$. Thus,
applying Corollary~\ref{cor:PQSurgeryFormulaSum} twice, one gets a
relationship between $\NInv(Z)$, $\NInv(Z')$ and the Alexander
polynomial for the zero-surgery. If the divisibility $d(X)$ of $X$
equals $n$, then we call the surgery from $Z$ to $Z'$ a surgery with
divisibility $n$.
\end{remark}

Let $L_{p,q}$ denote the lens space $S^3=\{(w,z)\in\C^2\big|
|w|^2+|z|^2=1\}$ modulo the equivalence relation $$(w,z)\sim
(e^{\frac{2\pi i}{p}}w,e^{\frac{2\pi q i}{p}}z).$$ In order to
calculate $\epsilon(p,q,1)$, we rely on a calculation given in
Section~\ref{sec:Lens}:

\begin{prop}
\label{prop:LensCalc}
$$2
\left(\sum_{\spinc\in\SpinC(L_{p,q})}\NInv_{L_{p,q}}(\spinc)\right)= -p
\cm s(q,p),$$ where $s(q,p)$ is the Dedekind sum 
$$s(q,p)=\frac{1}{4p}\sum_{k=1}^{p-1}
\cot\left(\frac{\pi k}{p}\right) \cot\left(\frac{\pi k q}{p}\right).$$
\end{prop}

As a consequence of Corollary~\ref{cor:PQSurgeryFormulaSum} and
Proposition~\ref{prop:LensCalc}, we have the following:

\begin{theorem}
\label{thm:CassonWalker}
$$2 \left(\sum_{\spinc\in\SpinC(Y)}\NInv_{Y}(\spinc)\right)=
\big|H_1(Y;\Z)\big|\lambda(Y),$$ where $\lambda(Y)$ is the
Casson-Walker invariant of $Y$.
\end{theorem}

\begin{proof}
Since any rational homology three-sphere $Y$ can be obtained from
$S^3$ by a sequence of surgeries on rational homology three-spheres,
Corollary~\ref{cor:PQSurgeryFormulaSum} and the constants $\epsilon(p,q,d)$
determine the sum
$$\sum_{\spinc\SpinC(Y)}\NInv_Y(\spinc)$$ 
(see Remark~\ref{rmk:Surgeries}). A suitably normalized version of 
the Casson-Walker invariant satisfies a formula of the same shape, with
constants $\epsilon'(p,q,d)$. We spell this out as follows.

Let $|\Tors(H_1(X))|=k$, so that $|\Tors H_1(Y_0)|=k/d$ and $|H_1(Y)|=dk$, 
$|H_1(Y_{p/q})|=pdk$.

Walker's surgery formula (see p.~82 of~\cite{Walker}) says that:
$$\lambda(Y_{p/q})=\lambda(Y)+\frac{q}{pd^2}\Gamma(X)-s(q,p)+\frac{(d^2-1)q}{12d^2p}.$$
Here, $\Gamma(X)$ is the sum $$\sum_j b_j j^2,$$ where $A_X(T)=\sum
b_j T^j$ is the symmetrized Alexander polynomial of $X$, normalized so
that $A_X(1)=1$. (The symmetrization forces us to allow half-integer
powers of $T$, if $d$ is even.)  We compare $\Gamma(X)$ with $\sum_j
a_j j^2$. By comparing with the Milnor torsion, one sees that $A_X =
A_{Y_0}(1+T+...T^{d-1})$ up to possible multiples of $T$ and
constants. Now, $a_j$ are the coefficients in the Alexander polynomial
of $Y_0$, normalized so that $A_{Y_0}(1)=|\Tors H_1(Y_0;\Z)|$. With
these normalizations, then, it follows that 
$$A_X = \frac{1}{k}A_{Y_0}(T)\left(\frac{T^{d/2}-T^{-d/2}}{T^{1/2}-T^{-1/2}}\right).$$
Then, 
\begin{eqnarray*}
\sum_j b_j j^2 &=& \frac{d^2}{dT^2}A_{X}(1) \\
&=& \frac{d}{k} (\frac{d^2}{dT^2}A_{Y_0}(1)) + \frac{1}{d} \sum_{\frac{-d+1}{2}}^{\frac{d-1}{2}} i^2 \\
&=& \frac{d}{k} \left(\sum_j a_j j^2\right) + 
\frac{(d^2-1)}{12}
\end{eqnarray*}

For the renormalized version $$\lambda'(Y)=\OneHalf
|H_1(Y)|\lambda(Y),$$ we then have: 
\begin{eqnarray*}\lambda'(Y_{p/q})&=& p
\lambda'(Y) + \frac{q}{2} \left(\sum_j a_j j^2\right) +
\frac{qk(d^2-1)}{12 d} - \frac{pkd\cm s(q,p)}{2} \\
&=& 
p \lambda'(Y) + q\left(\sum_{j\geq 1} a_j j^2\right) +
|\Tors H_1(X;\Z)| \left(\frac{q(d^2-1)}{12 d}-\frac{p d \cm s(q,p)}{2}\right).
\end{eqnarray*}

It follows that $\epsilon'(p,q,d)= \left(\frac{q(d^2-1)}{12 d}-\frac{p
d \cm s(q,p)}{2}\right)$.

Thus, it remains to show that
$$\epsilon(p,q,d)=\epsilon'(p,q,d).$$
For $d=1$ it follows from
Proposition~\ref{prop:LensCalc} that $2\epsilon(p,q,1)=-p\cm s(q,p)$, so we have that $\epsilon(p,q,1)=\epsilon'(p,q,1)$ as claimed.

We now argue that in fact $\epsilon(p,q,d)$ is determined by the
surgery formula and the values of $\epsilon(p,q,1)$. To this end, we
find it convenient to make the following definitions. Choose three
fiber circles in $S^2\times S^1$, and let $M(p_1,q_1,p_2,q_2,p_3,q_3)$
denote the manifold obtained by performing $p_i/q_i$ surgery on the
$i^{th}$ circle, with respect to the framing induced by the product
structure. Similarly, let $N(p_1,q_1,p_2,q_2)$ denote the manifold
obtained from surgeries along only two circles, and let
$Y(p_1,q_1,p_2,q_2)$ be the three-manifold obtained from
$N(p_1,q_1,p_2,q_2)$ by deleting a tubular neighborhood of the third
circle.
Note that $N(p_1,q_1,p_2,q_2)$ is either a lens space or $S^1\times
S^2$. 

Note that $M(n,1,-n,1,q,-p)$ is obtained from $M(n,1,-n,1,0,1)$ by a
$(p,q,n)$ surgery. Note also that $M(n,1,-n,1,0,1)$ is a connected sum of
lens spaces $L(n,1)\#{\overline{L(n,1)}}$. This is clear, for example,
from the Kirby calculus picture in Figure~\ref{fig:KirbyCalc}.


\begin{figure}
\mbox{\vbox{\epsfbox{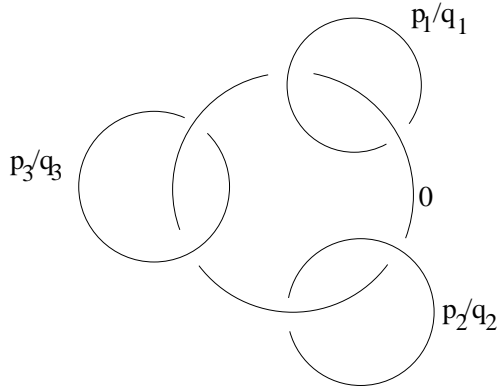}}}
\caption{\label{fig:KirbyCalc}
Kirby Calculus picture of $M(p_1,q_1,p_2,q_2,p_3,q_3)$}
\end{figure}

Suppose first that $q$ and $n$ are relatively prime. Then, since
$d(Y(n,1,q,-p))=1$, it follows that $M(n,1,-n,1,q,-p)$ is obtained
from $S^3$ by a sequence of surgeries of divisibility $1$.

With this in place, we turn our attention to the cases where $n=2$. We
can assume that $q$ is even. Consider the manifold $Y(2,1,2,1)$.
It has divisibility equal to $2$. It is easy to see that
$M(2,1,2,1,p-q,q)$ is gotten from $M(2,1,2,1,1,0)$ by a $(p,q,2)$
surgery.
Both of these manifolds can be obtained from $S^3$ by
surgeries of divisibility $1$. (Note that the first manifold is
$L(4,1)$, and $Y(2,1,p-q,q)$ has divisibility $1$, since $p-q$ is odd.)

For the general case, we use induction on $d$. Note that
$d(Y(n,1,q,-p))\leq n$, with equality iff $q=tn$ and
$n|t-p$. Similarly, $d(Y(-n,1,q,-p))\leq n$, with equality iff $q=tn$
and $n|t+p$. Now it follows that $M(n,1,-n,1,q,-p)$ can be obtained
from either $N(n,1,q,-p)$ or $N(-n,1,q,-p)$ with surgeries of
divisibility less than $n$, unless $q=tn$ and $n|2p$. Since $n>2$,
this would imply that $p$ and $q$ are not relatively prime.

Now, since $\epsilon(p,q,1)=\epsilon'(p,q,1)$, it follows from the
above argument that $\epsilon(p,q,d)=\epsilon'(p,q,d)$ for all
$d$. This finishes the proof of Theorem~\ref{thm:CassonWalker}.
\end{proof}

Theorem~\ref{thm:PQSurgeryFormula} rests on an analogue of
Proposition~\ref{prop:SurgerTheta}, which gives a relation among
$\theta_Y$, $\theta_{Y_{p/q}}$, $\theta_{Y_0}$, and certain spectral
flow terms. 
As before, we can define
$$\ModSp_{\met_t}(X)=\{(D,s,t)\in\Sym^{g-1}\times[0,1]\times[0,1]\big|
s\leq t, 
\Theta_{\met_s}(D)\in L(U_0),
\Theta_{\met_t}(D)\in\Lambda(U_1)\},$$
where the handlebodies $U_0$ and $U_1$ refer to the Heegaard
decomposition of $Y$. Note, however that the spaces $\Lag(U_0)$ and 
$\Lambda(U_1)$
depend only on $X$.

We describe how to partition this moduli space according to $\SpinC$
structures on $X$.  As in~\cite{Theta}, we consider the map
$$\begin{CD}
\pi_1(\Jac)\cong H^1(\Sigma;\Z)@>{\delta}>>H^2(X;\Z).
\end{CD}$$ Its
kernel gives rise to a covering space $\OurJac$ of $\Jac$, with
transformation group $H^2(X;\Z)$. A spin structure $\spinc$ on $X$
gives rise to lifts of $L(U_0)$ and $\Lambda(U_1)$ in $\OurJac$, up to
simultaneous translation by elements of $H^2(X;\Z)$ as follows. Let
$\spinc_0$ be a spin structure on $X$, and let $p\in \Jac$ be the
corresponding point. Any $\SpinC$ structure on $X$ can be written as
$\spinc_0+\ell$, where $\ell\in H^2(X;\Z)$. Let ${\widetilde p}$ be
any lift of $p$ to $\OurJac$. Then $L_0(\spinc)$ is the lift of $L_0$
to $\OurJac$ which passes through ${\widetilde p}$, and
$\Lambda_1(\spinc)$ is the lift of $\Lambda_1$ which passes through
${\widetilde p}+\ell$. (It is easy to see that these subspaces are
independent of the spin structure, as stated.) Note that the quotient
map, induces a diffeomorphism of each $L_0(\spinc)$,
resp. $\Lambda_1(\spinc)$, to the corresponding spaces $L_0$ and
$\Lambda_1$ respectively. Also, there is a lifting
$$\OurTheta_\met\colon
\OurSym \longrightarrow \OurJac,$$ where $\OurSym$ is the lift of
$\Sym^{g-1}(\Sigma)$ corresponding to kernel of the map
$$\begin{CD}
\pi_1(\Sym^{g-1}(\Sigma))@>>>
H_1(\Sym^{g-1}(\Sigma))@>{(\Theta_\met)_*}>>
H_1(\Jac)\cong H^1(\Sigma;\Z) @>{\delta}>> H^2(X;\Z).
\end{CD}
$$

Correspondingly, let
$$\ModSp_{\met_t}(X,\spinc)
=\{(D,s,t)\in\Sym^{g-1}\times[0,1]\times[0,1]\big|
s\leq t, 
\OurTheta_{\met_s}(D)\in L_0(\spinc),
\OurTheta_{\met_t}(D)\in\Lambda_1(\spinc)\}.$$
These spaces naturally give a partition
$$\ModSp_{\met_t}(X)=\coprod_{\spinc\in\SpinC(X)}\ModSp_{\met_t}(X,\spinc).$$

As in Section~\ref{sec:Casson}, we have a map: $$\rho\colon
\ModSp_{\met_t}(X)\longrightarrow
\Jac(\partial X)\cong \Torus{2},$$
defined by
$$\rho(D,s,t)=\Hol_{m\times\ell}(\Theta_{\met_t}(D)),$$ normalized so
that the point $(\OneHalf,\OneHalf)\in\Torus{2}$ corresponds to the
spin structure on $S^1\times S^1$ which does not bound, and the 
circle of reducibles
for any  level zero $\SpinC$ structure 
on $X$ restrict to give the circle $S^1\times\{0\}$.

\begin{prop}
\label{prop:StructureModSpaceQ}
For any $\epsilon>0$, if $\met_1$ is sufficiently stretched out normal
to the attaching circles $\{\beta_1,...,\beta_{g-1}\}$, then the
moduli space $\ModSp_{\met_t}(X,\spinc)$ is generically a compact,
smooth, one-dimensional manifold with two types of boundary components
corresponding to $t=1$ and $s=t$.  Furthermore, the $t=1$ boundary
maps under $\rho$ into an $\epsilon$-neighborhood of
$\OneHalf\times\OneHalf$, and those with $s=t$ map under $\rho$ to
the circle
$\gamma_{0}(y)=S^1\times \{\frac{y}{d}\}$. 
\end{prop}

\begin{proof}
The $s=t$ boundaries correspond to the intersection of the theta divisor with
the $L(U_0)\cap \Lambda(U_1)$. Using the identification between the
Jacobian and the $H^1(\Sigma;S^1)$, $L(U_0)\cap \Lambda(U_1)$
corresponds to those $S^1$ representations of $\pi_1(\Sigma)$ which
extend to representations of $\pi_1(X)$. Since $d\ell$ bounds in $X$,
these representations must take $\ell$ to a $d$-torsion
point. Tensoring with any element of $H^2(Y_0;\Z)$ which maps to a
generator of $H^2(Y_0;\Z)/\Tors$, and taking the reducible representative
on $X$, changes the holonomy around $\ell$ by $1/d$: this shows that
the number $y$ appearing above is the level of the $\SpinC$ structure
as defined in the beginning of this section.
\end{proof}

Let $\gamma_{p/q}$ denote the circle in $S^1\times S^1$ with slope
$p/q$ and which goes through the point $(0,0)$ if $q$ even and
$(\OneHalf,0)$ if $q$ is odd. Let $\gamma_0(y)$ denote the circle
$S^1\times \{\frac{y}{d}\}$, and $\gamma_{0}(y)'$ denote the 
vertical translate of $\gamma_{0}(y)$ by some small
$\delta>0$. Let $\gamma$ be the
circle $0\times S^1$.

The curve $\gamma_{p/q}$ meets $\gamma_0(y)$ in $p$ points. Using the
orientation on $H^1(X;\R)$ coming from $m$, we can label these in
increasing order $\phi_1<\phi_2<...<\phi_p$. Also, $\gamma_0(y)$
intersects $\gamma$ in a single point $\phi_0=(0,y)$. To each $\phi_i$,
we can associate the interval $[\phi_i,\phi_0]$. Now, to any interval
$I\subset \gamma_0(y)$, we can associate a spectral flow
$\SF_{\R\times \Sigma}I$ defined as follows.  Each interval is covered
by $d$ one-parameter families of (nearly reducible)
$\SpinC$ connections with traceless
curvature in the $\SpinC$ structure $\spinc$ on $X$. By restricting to
$\R\times\Sigma\subset X$, we get $d$ one-parameter families of Dirac
operators. $\SF_{\R\times
\Sigma}I$, then, denotes the sum of these spectral flows.

We have the following analogue of Proposition~\ref{prop:SurgerTheta}.
To state it, we must make some preliminary remarks concerning spectral
flow on $X$ and $\R\times \Sigma$. The set of $\SpinC$ connections on
$\spinc$ over $X$ with traceless curvature is parameterized by the
circle $\frac{H^1(X;\R)}{H^1(X;\Z)}$, which we orient via the homology
class $m\in H_1(X;\Z)$. These operators give rise to a family
$\{A_t\}$ of Dirac operators over $\R\times\Sigma$. 
Strictly speaking,
this family might not be a family of Fredholm operators, if the level
is $d/2$, since in this case the boundary value of $\{A_t\}$ goes
through the bad point. We can compensate for this by introducing a
some curvature so that the the boundary value maps to a curve of the
form $S^1\times (\OneHalf+\epsilon)$.

\begin{lemma}
\label{lemma:SFCircle}
The spectral flow around the circle determined by $\{A_t\}$
depends only $y$ and
$d$. 
\end{lemma}

\begin{proof}
This follows immediately from the
splitting principle for spectral flow, applied to the zero-surgery
$Y_0$.
\end{proof}

\begin{prop}
\label{prop:SurgerThetaPQ}
Fix a generic path of metrics $\met_t$ so that
$\rho(\ModSp_{\met_t}(X,\spinc))$ misses the points $\phi_0,\phi_1,...,\phi_p$.
Then, we have that
\begin{eqnarray}
\sum_{\spincb\in\SpinC(Y_{p/q};\spinc)}
\Inv_{Y_{p/q}}(\spincb)&=&
p \left(\sum_{\spincc\in\SpinC(Y;\spinc)}\Inv_{Y}(\spincc)\right)
+ q
\left(\sum_{\spincd\in\SpinC(Y_0;\spinc)}\Inv_{Y_0}(\spincd)\right) 
\nonumber \\
&&
+ \sum_{i=1}^p \SF_{\R\times\Sigma}[\phi_i,\phi_0]
+f(p,q,d,y),
\label{eq:SurgerThetaPQ}
\end{eqnarray}
where $d=d(X)$ and $y=y(\spinc)$, and $f(p,q,d,y)$ is an integer which
is independent of $X$.
\end{prop}

\begin{proof}
Note first of all that $\gamma_{p/q}$ is homologous to
$p\gamma+q\gamma_0'(y)$. Thus, we can
consider the oriented $2$-chain $\Dom \subset S^1\times S^1$ which does not
contain $\OneHalf\times \OneHalf$, and whose boundary is 
$\gamma_{p/q}-p\gamma-q\gamma_0'(y)$. By elementary differential topology,
$$\partial (\Restrict^{-1}(\Dom))=\Restrict^{-1}(\gamma_{p/q})-
p \Restrict^{-1}(\gamma)- q \Restrict^{-1}(\gamma_0'(y)) + 
\left(\Restrict|_{\partial\ModFlow_{\met_t}(X;\spinc)}\right)^{-1}(\Dom).$$
These are to be viewed $0$-chains; i.e. 
we have that 
\begin{equation}
\label{eq:PreSurgery}
\#\left(\partial (\Restrict^{-1}(\Dom))=\#\Restrict^{-1}(\gamma_{p/q})\right)-
p \left(\#\Restrict^{-1}(\gamma)\right)- q \left(\#\Restrict^{-1}(\gamma_0'(y))\right) + 
\sum_{x\in \partial\ModFlow_{\met_t}(X;\spinc)}
u_x
\end{equation}
where $u_x$ denotes the multiplicity of $\Dom$ at $\Restrict(x)$
times the degree of $\Restrict$ to $x$.  Note that
\begin{eqnarray}
\#\Restrict^{-1}\gamma_{p/q} &=& 
\sum_{\spincb\in\SpinC(Y_{p/q};\spinc)}
\Inv_{Y_{p/q}}(\spincb)
\label{eq:SurgYPQ}\\
\#\Restrict^{-1}(\gamma)&=&
\sum_{\spincc\in \SpinC(Y;\spinc)} \Inv_{Y}\left(\spincc\right) 
\label{eq:SurgY}
\\
\#\Restrict^{-1}(\gamma_0'(y)) &=&
\sum_{\spincd\in\SpinC(Y_0;\spinc)}\Inv_{Y_0}(\spincd).
\label{eq:SurgYZ}
\end{eqnarray}

Note that the multiplicity of $\Dom$ is constant on the open interval
$(\phi_i,\phi_{i+1})$ and $(\phi_p,\phi_0)$, and we denote the constant by
$\nu_i$ and $\nu_p$ respectively. 
Thus, we can write the final term in
Equation~\eqref{eq:PreSurgery} as:
\begin{eqnarray*}
\sum_{x\in \partial\ModFlow_{\met_t}(Y-K;\spinc)}
u_x &=&
\left(\sum_{i=0}^{p-1} \nu_i~\SF_{\Cyl}[\phi_i,\phi_{i+1}]\right)
+ \nu_p~\SF_{\Cyl}[\phi_p,\phi_0],
\end{eqnarray*}
thanks to Proposition~\ref{prop:SFCylBig}. Moreover, it is a simple
homological fact that $\nu_i=\nu_0 + i$, so we can rewrite the above
quantity as a sum of 
$$\sum_{i=1}^p \SF_{\Cyl}[\phi_i,\phi_0]
+ \nu_0 f(\spinc),$$
where $f(\spinc)$ is the spectral
flow around the circle from Lemma~\ref{lemma:SFCircle}.
Substituting this back into
Equation~\eqref{eq:PreSurgery}, along
Equations~\eqref{eq:SurgY}-\eqref{eq:SurgYZ}, we obtain
Equation~\eqref{eq:SurgerThetaPQ}. 

Note that in the case where the level is $2$, the spectral flow terms 
might not be defined, since the path of reducibles fails to induce a 
Fredholm operator. This can be compensated by using spectral flow 
through nearly-reducible connections (i.e. translating $\gamma_{0}$, 
once again, by a small amount).
\end{proof}

The next step is to see how $\Corr$ changes under surgeries, which is
another application of the excision principle for indices. Fix
$\SpinC$ structures $\spincc$ and $\spincb$ on $Y_{p/q}$ and $Y$
which restrict to the same $\SpinC$ structure on $Y-K$ (so that the
corresponding $\SpinC$ connections induced on $X$ can be connected through
connections with traceless curvature). 

There is a standard cobordism $Z_{p/q}$ between the neighborhood of a
knot, thought of as the core of $U=\CDisk\times S^1$, and its $p/q$
surgery, $U_{p/q}$. Of course, $U_{p/q}$ is also diffeomorphic to $U$,
but we will use different metrics. From standard Kirby calculus, this
cobordism is obtained as a plumbing of two-spheres, prescribed by the
Hirzebruch-Jung fraction expansion of $p/q$. Since the original
manifolds have boundaries, the ``cobordism'' $Z_{p/q}$ is actually a
manifold-with-corners, with a boundary $[0,1]\times S^1\times S^1$,
and two others which are $U$ and $U_{p/q}$.  In particular, if we fix
$S^1\times \CDisk$, then $\left([0,1]\times S^1\times
\CDisk\right)\cup_{[0,1]\times S^1\times S^1} Z_{p/q}$ is the standard
cobordism between $S^3$ and the lens space $L_{p,q}$ and, more
generally, if $X=Y-K$, then the manifold
$$W_{p/q}=\left([0,1]\times X\right)\times_{[0,1]\times S^1\times S^1} Z_{p/q}$$ is
a cobordism between the manifold $Y$ and $Y_{p/q}$, identifying
$S^1\times S^1$ with the meridian times the longitude.

We endow $Z_{p/q}$ with a metric compatible with its
manifold-with-corners structure; i.e. so that it is product-like in
neighborhoods of each of its boundaries. Indeed, we can choose the
metric on $[0,1]\times S^1\times S^1$ to be independent of the $[0,1]$
factor. If $X$ is fixed with a product-like metric near its boundary,
where it is identified with $S^1\times S^1$ (meridian times
longitude), then this, together with a neck-length parameter $T$,
gives both $Y$ and $Y_{p/q}$ preferred metrics.

Fix a $\SpinC$ structure $\spincr$ on $W_{p/q}$ with
$\spincr|{Y_{p/q}}=\spincb$, and $\spincr|{Y}=\spincc$. We can choose
a $\SpinC$ connection $A$ on $\spincr$,
whose restriction to $[0,1]\times (Y-K)$ induces a one-parameter
family of reducibles $\{A_t\}$ on $Y-K$.
We then have the following analogue of Lemma~\ref{lemma:SurgerCorr}:

\begin{prop}
\label{prop:SurgerCorrPQ}
Let $\spinca=\spincb|(Y-K)=\spincc|(Y-K)$. 
The difference $\Corr(Y_{p/q},\spincb)-\Corr(Y,\spincc) - 
\SF_{Y-K}(\{A_{t}\})$ is independent of the
manifold $Y$, depending only on $p$, $q$, $d$, and $y(\spinc)$ 
\end{prop}

\begin{proof}
Since $\SF_{Y-K}(A_t)$ can be thought of as the $L^2$
index of the Dirac operator for $(Y-K)\times [0,1]$, with ends
attached, the above statement is an easy application of the splitting
principle for the index, bearing in mind that the difference
$\Corr(Y_{p/q},\spincb)-\Corr(Y,\spincc)$ is the index of the Dirac
operator on $W_{p/q}$, and topological terms which depend only on
$p$, $q$, $d$, and $y(\spinc)$. 
\end{proof}

\vskip.2cm
\noindent{\bf{Proof of Theorem~\ref{thm:PQSurgeryFormula}.}}
The surgery formula for $\theta$ shows that the the invariant
transforms in the manner which $\NInv$ is claimed to transform under,
plus a sum of spectral flow terms on $\R\times \Sigma$. We show (with
the help of Proposition~\ref{prop:SurgerCorrPQ}) that these spectral
flow terms are cancelled by the correction terms $\CorrU$ for the
handlebodies and the terms for $\CorrY$ for the manifolds $Y$ and
$Y_{p/q}$. This is true since first of all the 
$$\CorrU(U_1(p/q),\spincb)-\CorrU(U_1,\spincc) =
\SF_{U_1-K}(\spincc,\spincb) + \left(\CorrU(U(p/q))-\CorrU(U)\right),$$
so
\begin{eqnarray*}
\lefteqn{\left(\CorrU(U_0,\spincb) + \CorrU(U_1(p/q),\spincb)\right)
- \left(\CorrU(U_0,\spincc) + \CorrU(U_1,\spincc)\right)} \\
&=& \SF_{U_0}(\spincc,\spincb) + \SF_{U_1-K}(\spincc,\spincb) +
\left(\CorrU(U(p/q))-\CorrU(U)\right)
\end{eqnarray*}
Applying the splitting theorem for spectral flow
$$\SF_{Y-K}(\spincc,\spincb)=\SF_{U_0}(\spincc,\spincb)+
\SF_{\R\times\Sigma}(\spincc,\spincb)
+ \SF_{U_1-K}(\spincc,\spincb),$$
Proposition~\ref{prop:SurgerThetaPQ} and
Proposition~\ref{prop:SurgerCorrPQ}, the theorem follows.
\qed


\section{The invariant $\NInv$ for lens spaces}
\label{sec:Lens}

In this section, we calculate the $\NInv$-invariant for lens spaces
$L_{p,q}$. Our aim is the following:

\begin{prop}
\label{prop:LensSpaces}
There is an identification
$\SpinC(L_{p,q})\cong \Zmod{p}$, under which
the invariant $\NInv$ for $L(p,q)$ is given by
$$2~\NInv_{L_{p,q}}(\alpha)=-s(q,p)-\frac{1}{2p}\sum_{g=1}^{p-1}
\csc\left(\frac{\pi g}{p}\right) \csc\left(\frac{\pi q g}{p}\right)
\cos\left(\frac{2\pi g \alpha}{p}\right)$$
if $p$ is even and 
$$2~\NInv_{L_{p,q}}(\alpha)=-s(q,p)-\frac{1}{2p}\sum_{g=1}^{p-1}
\csc\left(\frac{2 \pi g}{p}\right) \csc\left(\frac{2 \pi q g}{p}\right)
\cos\left(\frac{2\pi g \alpha}{p}\right)$$
if $p$ is odd.
In particular, in either case,
\begin{equation}
\label{eq:TotalInvLens}
2\sum_{\spinc\in \SpinC(L_{p,q})}
\NInv_{L_{p,q}}(\spinc) = -p\cdot s(q,p).
\end{equation}
\end{prop}

\begin{remark}
The reader can find the precise identification $\SpinC(L_{p,q})\cong
\Zmod{p}$ in the course of the proof of the above result (see the
proof of Lemma~\ref{lemma:EtaDirac}). We do not spell this out at this
point, since our main interest is Equation~\eqref{eq:TotalInvLens}
which is independent of these identifications.
\end{remark}

Note first that $L_{p,q}$ can be thought of as a quotient of the
standard, round three-sphere $S^3$ by a group of isometries. We call
this the standard metric on $L_{p,q}$. We can work with this metric to
calculate $\NInv$, according to the following lemma, whose proof is
given at the end of this section:

\begin{lemma} 
\label{lemma:HeegaardMetric}
On a lens space $L_{p,q}$, we have that
$$\NInv(L_{p,q},\spinc)=-\CorrY(L_{p,q},\spinc)$$ for any $\SpinC$
structure $\spinc\in\SpinC(L_{p,q})$, where the right hand side is
calculated using the standard metric on $L_{p,q}$.
\end{lemma}

The calculation of $\Corr(L_{p,q})$ is a straightforward consequence
of the Atiyah-Bott-Lefschetz fixed point theorem, for manifolds with
boundary (a variant of which is mentioned in~\cite{APSII}, with a
complete proof given in~\cite{Donnelly}). Closely related calculations
can be found in (\cite{Gilkey}, see also~\cite{RationalBlowdown}). We
review the relevant theory.

We find it convenient to express the correction terms $\CorrY(Y)$ in
terms of Atiyah-Patodi-Singer eta-invariants, by the APS index
theorem. If $(X,\spincX)$ is a $\SpinC$ four-manifold which bounds a
rational homology three-sphere $Y$ with $\SpinC$ structure $\spinc$,
$$\Ind_{\R}(\Dirac,\spincX)= - \frac{1}{12}\int_X p_1(X) +
\frac{1}{4}\int_X c_1(\spincX)^2 - \frac{\etaDirac(Y,\spinc)+h}{2},$$
(where $h$ is the real dimension of the kernel of the Dirac operator on
$Y$, and $\etaDirac(Y)$ is its real eta invariant), 
while $$\sgn(X)= \frac{1}{3}\int_X p_1(X) -\etaSign(Y).$$ Here,
$p_1(X)$ is the first Pontryagin form of $X$ and $\etaSign(Y)$ is the eta
invariant for the signature operator for even forms, appearing in
Theorem~4.14 of~\cite{APSI}.  (Metrics are suppressed from the
notation, but one should bear in mind that the index of
the Dirac operator, the first Pontryagin form,
and the eta invariants and $h$ of the boundary all depend on metrics;
and moreover the metric on $X$ must have a cylindrical collar near its boundary
for the formula to hold.)  It follows, then that, if the Dirac
operator on $Y$ has no kernel, then
\begin{equation}
\label{eq:ReexpressCorrY}
\CorrY(Y)=\OneHalf \Ind_{\R}(\Dirac,\spincX) -
\frac{\langle c_1(\spincX)^2, [X]\rangle}{8} + \frac{\sgn(X)}{8} =
-\frac{1}{4}\etaDirac(Y,\spinc) - \frac{1}{8}\etaSign.
\end{equation}

For the standard metric on $L_{p,q}$, the formula 
\begin{equation}
\label{eq:EtaSignature}
\etaSign(L_{p,q})=
-4s(q,p)
\end{equation}
can be found in Proposition~2.12 of~\cite{APSII}. This
formula, and the formula for $\etaDirac$, are both obtained by a
Lefschetz-type theorem, which we outline for Dirac operator.

Let $g$ be an isometry of ${\widehat Y}$, which is lifted to an
automorphism of the spinor bundle of ${\widehat Y}$. Then, if
$V_\lambda$ denotes the $\lambda$-eigenspace of the Dirac operator,
$g$ induces an automorphism of $V_\lambda$, and there is an associated
$\eta$-function $\eta_g(s,{\widehat Y})$, defined by analytically
continuing the function of $s$ $$\eta_g(s,{\widehat
Y})=\sum_{\lambda\in (\Spec\Dirac-0)}\frac{({\mathrm{sign}} \lambda)
\Tr (g|V_\lambda)}{|\lambda|^s}$$ over the complex plane.  
The $g$-eta invariant, then, is defined to be the evaluation
$\eta_g({\widehat Y})=\eta_g(0,{\widehat Y})$.
Moreover, if $G$ is a finite group of automorphisms of a spin-manifold
${\widehat Y}$, which acts freely by isometries on the base, and
$\alpha$ is any representation of $G$, there is an associated flat
line bundle $F_{\alpha}$ over the spin manifold $Y={\widehat
Y}/G$. Coupling the Dirac operator to $F_{\alpha}$, we obtain an
$\eta$-invariant $\eta_{\alpha}(Y)$ obtained by evaluating
$\eta_{\alpha}(s,Y)$ at $s=0$. If the representation maps to $S^1$,
this is the eta invariant for the Dirac operator on $Y$ coupled to the
$\SpinC$ structure obtained by tensoring the spin bundle on $Y$ with
$F_\alpha$. The eta invariants are related by:
\begin{equation}
\label{eq:RelateEtas}
\eta_\alpha(s,Y)=\frac{1}{|G|}\sum_{g\in G}\eta_g(s,{\widehat Y})
\chi_{\alpha}(g),
\end{equation}
where $\chi_{\alpha}$ is the character of $\alpha$. 

If the action of $g$ on ${\widehat Y}$ extends to $X$ as an isometry acting on
the spinor bundle, then there is a $g$-index $$\ind(D,g)=\Tr(g|\Ker
D)-\Tr(g|\Coker D).$$ The APS version of the $G$-index theorem (see
Theorem~1.2 of~\cite{Donnelly}) gives a formula for the $g$-index in
terms of local numbers around the fixed points of $g$, and the $g$-eta
invariant. If the fixed point set of $g$ on $X$ is isolated, it takes
the particularly simple form:
\begin{equation}
\label{eq:GIndexThm}
\ind(D,g)= \sum_{x\in {\mathrm Fix}(g)}
\frac{\Tr(g|W^+_x)-\Tr(g|W^-_x)}{\det(1-d_xg)} - \frac{\eta_g(0)+h_g}{2},
\end{equation}
where $d_xg$ is the differential of $g$ at $x$.

\begin{lemma}
\label{lemma:EtaDirac}
For the standard metric on $L_{p,q}$, the eta invariant for the Dirac
operator is given by
$$\etaDirac_{\alpha}(0)=
\left\{\begin{array}{ll}
-\frac{1}{p}\sum_{g=1}^{p-1}
\csc\left(\frac{\pi g}{p}\right) \csc\left(\frac{\pi q g}{p}\right)
\cos\left(\frac{2\pi g \alpha}{p}\right) & {\text{if $p$ is even}} \\
-\frac{1}{p}\sum_{g=1}^{p-1}
\csc\left(\frac{2\pi g}{p}\right) \csc\left(\frac{2\pi q g}{p}\right)
\cos\left(\frac{2\pi g \alpha}{p}\right) & {\text{if $p$ is odd}}
\end{array}
\right..
$$
\end{lemma}

\begin{proof}
The $\eta$-invariants of $L_{p,q}$ can be calculated using $G$-index
theorem, applied to a $\Zmod{p}$ action on the four-ball $B^4$ endowed
with an $SO(4)$-invariant metric with non-negative scalar curvature
which is product-like near the boundary. To lift the action to the
spinor bundle, we follow~\cite{Austin}. Let $h$ be the rotation
$$h(w,z)=(\zeta w,
\zeta^q z),$$ where $\zeta_p$ is a primitive $p^{th}$ root of
unity. The element $h$ generates a $\Zmod{p}$ action on $B^4$. This
can be lifted to a map ${\widetilde h}$ which acts on the spinor
bundles $W^+$ and $W^-$ by choosing a square root $\gamma$ of $\zeta$
and letting $${\widetilde g}(\Phi_\pm )= \gamma^{\pm 1-q} \Phi_\pm,$$
where $\Phi_\pm$ is a spinor in $W^\pm$, viewed as a bundle of
quaternions. There are two slightly different cases, according to the
parity of $p$. If $p$ is even, then $\pm 1-q$ is even, so $\gamma^{\pm
1-q}$ has order $p$; if $p$ is odd then $(\gamma^{\pm 1 - q})^2$ has
order $p$. Thus, letting ${\widetilde g}={\widetilde h}$ if $p$ is
even and ${\widetilde g}={\widetilde h}^2$ if $p$ is odd, we see that
multiples of ${\widetilde g}$ generate a free $\Zmod{p}$ action on
$B^4$ together with its spinor bundle.

We apply the $G$-index theorem to any $g\in \Zmod{p}$, to
compute the $g$-eta invariant. First note that the
four-ball has non-negative scalar curvature, so the kernel and
cokernel of the Dirac operator vanish; in particular, the $g$-index
vanishes. 
If $g\neq 0$, the only fixed point of $g$ is the
origin, and it is easy to calculate the contribution of that fixed
point to be 
\begin{eqnarray*}
-\frac{1}{2}\csc\left(\frac{\pi
g}{p}\right)\csc\left(\frac{\pi q g}{p}\right)
&{\text{resp.}}& 
-\frac{1}{2}\csc\left(\frac{2\pi
g}{p}\right)\csc\left(\frac{2\pi q g}{p}\right)
\end{eqnarray*}
according to whether $p$ is even or odd. It then follows
from the $G$-index theorem that, for $g\neq 0$,
$$\eta_g(0)=-\csc\left(\frac{\pi
g}{p}\right)\csc\left(\frac{\pi q g}{p}\right).$$ From
Equation~\eqref{eq:RelateEtas}, it follows that 
$$\eta_{\alpha}(0)=
-\frac{2}{p}\left(\eta(0,S^3)+ 
\sum_{g=1}^{p-1} \eta_g(0)
\cos\left(\frac{2\pi g \alpha}{p}\right)\right).$$
Since $S^3$ has symmetric spectrum, the term $\eta(0,S^3)=0$.
The lemma follows.
\end{proof}

We now give the proof of Lemma~\ref{lemma:HeegaardMetric}, to justify
our use of the ``standard metric''.

\vskip.2cm
\noindent{\bf{Proof of Lemma~\ref{lemma:HeegaardMetric}}}.
Suppose we have a genus one Heegaard decomposition
$L_{p,q}=U_0\cup_{S^1\times S^1} U_1$, and fix a pair $\Umet_0'$ and
$\Umet_1'$ of metrics on $U_0$ and $U_1$ with non-negative scalar
curvature metrics which bound a fixed flat metric on the torus. These
induce a metric $\Umet_0'\#_T \Umet_1'$ on $L_{p,q}$ with non-negative
scalar curvatures.  Then, for all $T$,
$$\NInv(\spinc)=
\CorrU(\Umet_0')+\CorrU(\Umet_1')-\CorrY(\spinc,\Umet_0'\#_T\Umet_1').$$
This is true since all $\Inv(\spinc)$ vanishes identically, as the
theta divisor on a torus does not bound.  Moreover,
$\CorrY(\spinc,\Umet_0'\#_T\Umet_1')$ is independent of $T$, since for
all $T$, the metrics have non-negative scalar curvatures, so the Dirac
operator never acquires kernel.

We connect the standard metric on $L_{p,q}$ with a metric of the form
$\Umet_0'\#_T\Umet_1'$ through a path of metrics with non-negative
scalar curvature, to show that the correction terms $\CorrY$ agree.
In fact, consider the one-parameter family of metrics on $S^3$
constructed in the proof of Proposition~\ref{prop:CompareCorrections}.
These metrics are always invariant under the $S^1\times S^1$ action on
$S^3$ (rotating the $\theta$ and $\phi$ coordinates). In particular,
they are invariant under the $\Zmod{p}\subset S^1\times S^1$ action
whose quotient gives $L_{p,q}$.
The requisite metrics on $L_{p,q}$, then are the metrics induced on
the quotient by the $\Zmod{p}$ action. 

Note that $\CorrU(\Umet_0')=\CorrU\left((\CDisk\times
S^1)/(\Zmod{p})\right)=0$, since the metric on $\CDisk\times
S^1/(\Zmod{p})$ can be realized metrically as a fiber bundle over
$S^1$ with fiber a disk (endowed with a circle-invariant metric with
non-negative sectional curvature), where the holonomy is rotation
through some angle. Since any rotation is isotopic to the identity
through isometries, we can connect $\CDisk\times S^1/(\Zmod{p})$ with
a product metric $\CDisk\times S^1$ through metrics of non-negative
sectional curvature. Hence, $\CorrU\left((\CDisk\times
S^1)/(\Zmod{p})\right)=0$. Similarly, $\CorrU(\Umet_1')=0$.  Thus, the
lemma follows.
\qed

Thus, we have given all the ingredients to
Proposition~\ref{prop:LensSpaces}. Specifically, we have:

\vskip.2cm
\noindent{\bf{Proof of Proposition~\ref{prop:LensSpaces}}.}
The proposition follows from Lemma~\ref{lemma:HeegaardMetric},
Equations~\eqref{eq:ReexpressCorrY} and \eqref{eq:EtaSignature}, and
then Lemma~\ref{lemma:EtaDirac}.
\qed

\commentable{
\bibliographystyle{plain}
\bibliography{biblio}
}

\end{document}